\documentclass[11pt,a4]{article}

\usepackage{graphicx}%
\usepackage{multirow}%
\usepackage{amsmath,amssymb,amsfonts}%
\usepackage{amsthm}%
\usepackage{mathrsfs}%
\usepackage[title]{appendix}%
\usepackage{xcolor}%
\usepackage{textcomp}%
\usepackage{manyfoot}%
\usepackage{booktabs}%
\usepackage{algorithm}%
\usepackage{algorithmicx}%
\usepackage{algpseudocode}%
\usepackage{listings}%
\usepackage{plain}
\usepackage{natbib}
\usepackage[colorlinks=true,linkcolor=red,citecolor=blue]{hyperref}

\theoremstyle{thmstyleone}%
\theoremstyle{thmstyletwo}%

\theoremstyle{thmstylethree}%

\begin{document}

\title{Influence of gauges in the numerical simulation of the time-dependent Ginzburg-Landau model}


\author{
	Cyril Tain$^{1,2}$, Jean-Guy Caputo$^2$, and Ionut Danaila$^{1, *}$ 
	\\ \\
	$^1$Univ Rouen Normandie, CNRS,  \\
	Laboratoire de Math{\'e}matiques Rapha{\"e}l Salem,\\
	UMR 6085, F-76000 Rouen, France\\
	$^2$Laboratoire de Math{\'e}matiques de l'INSA Rouen Normandie, UR 322,\\
	 685 Av. de l'Universit{\'e}, Saint-{\'E}tienne-du-Rouvray, 76800 France,\\
	$^*$ Corresponding author: ionut.danaila@univ-rouen.fr
}

\date{\today}
\maketitle

%
%


\abstract{The time-dependent Ginzburg-Landau (TDGL) model requires the choice of a gauge for the problem to be mathematically well-posed. In the literature, three gauges are commonly used: the Coulomb gauge, the Lorenz gauge and the temporal gauge. It has been noticed [J. Fleckinger-Pell{\'e} \textit{et al.}, Technical report, Argonne National Lab. (1997)] that these gauges can be continuously related by a single parameter considering the more general $\omega$-gauge, where $\omega$ is a non-negative real parameter. In this article, we study the influence of the gauge parameter $\omega$ on the convergence of numerical simulations of the TDGL model using finite element schemes. A classical benchmark is first analysed for different values of $\omega$ and artefacts are observed for lower values of $\omega$. Then, we relate these observations with a systematic study of convergence orders in the unified $\omega$-gauge framework.  In particular, we show the existence of a tipping point value for $\omega$, separating optimal convergence behaviour and a degenerate one. We find that numerical artefacts are correlated to the degeneracy of the convergence order of the method and we suggest strategies to avoid such undesirable effects. New 3D configurations are also investigated (the sphere with or without geometrical defect).}




\maketitle


\section{Introduction}

The time-dependent Ginzburg-Landau model (TDGL) is used to describe the dynamics of vortices in superconductors. It has been mathematically studied since the 1990's. In \cite{du1994} the definitions of gauges are introduced to ensure a mathematically well-posed problem. 
Three main gauges are commonly used to relate the main variables of the model: the electric potential $\phi$, the magnetic vector potential $\mathbf{A}$ and the quantum complex order parameter $\psi$.  The Lorenz gauge states that $\phi = -\operatorname{div}\mathbf{A}$, the Coulomb gauge  that $\operatorname{div}\mathbf{A} = 0$ and the temporal gauge that $\phi = 0$.

For a given gauge, discretizations of the TDGL equations have been investigated with different methods and a large volume of literature exists in this field (for a review of these studies see \cite{du2005} and the references therein). In particular, finite element (FE) methods have been extensively studied for the three gauges:
\begin{itemize}
	\item The Lorenz gauge was studied in \cite{gaoSun2014} using a Crank-Nicolson scheme in time and Lagrange FE in space. For squared geometries, the authors  observed singularities  for the magnetic field at corners. To avoid such singularities (which are numerical artefacts), a mixed scheme for the Lorenz gauge was suggested  in \cite{gaoSun2015}. By introducing $\operatorname{\textbf{curl}} \mathbf{A}$ as a supplementary unknown, the authors showed that the magnetic field was computed correctly and numerical artefacts were avoided. The convergence of the mixed scheme was proved in \cite{gaoSun2018} in general domains, including two-dimensional non-convex polygons. 
	\item Existence and uniqueness for the TDGL under the Coulomb gauge was studied in \cite{tang1995time}. A numerical analysis for the Coulomb gauge can be found in \cite{gao2023finite}. The authors used a backward Euler scheme in time. The vector potential $\mathbf{A}$ was approximated by lowest order Nedelec FE, $\phi$ and $\psi$ by linear Lagrange FE.
	\item For the temporal gauge, a backward Euler scheme in time and piecewise quadratic finite elements in space were used to solve the TDGL equations in  the pioneering work of \cite{du1994b}. A drawback of the temporal gauge when compared to the Lorenz gauge is the degeneracy of the parabolic equation for the vector potential $\mathbf{A}$. As a result, the convergence for $\mathbf{A}$ is one order lower than in the Lorenz gauge. 
\end{itemize}

To assess on the best adequacy of a finite-element scheme with the choice of the gauge, a comparison between three numerical methods was presented in 
\cite{gao2016}: two schemes written with Lagrange FE (one in the Lorenz gauge, the other in the temporal gauge) and a third scheme with a mixed formulation (with also Lagrange FE) using the temporal gauge. When  the temporal gauge was used, their results showed the degeneracy of the convergence order for the vector potential and its divergence. 

In this article, we address the question of the  selection of the optimal gauge for a finite-element discretization by studying 
the more general $\omega$-gauge defined as $\phi = -\omega\operatorname{div}\mathbf{A}$. This gauge was theoretically introduced in \cite{fleckinger96}, but not analysed numerically. It allows one to link the temporal gauge ($\omega = 0$) and the Lorenz gauge  ($\omega = 1$) continuously. We first  estimate the dependence of convergence rates on the value of  $\omega$ using specially designed manufactured solutions in two (2D) or three (3D) dimensions of space. Different types of finite-element discretizations are tested: Lagrange and Raviart-Thomas mixed FE  schemes for 2D, Lagrange and Raviart-Thomas-Nedelec mixed FE schemes for 3D. We then apply the generalized  $\omega$-gauge  to well known benchmarks for the TDGL problem and point out that numerical artefacts observed in some simulations are related to the degeneracy of convergence orders.  This study thus offers a unified framework that directly and continuously relate the influence of the gauge to the convergence of the FE numerical scheme.

The outline of the paper is as follows. In Sec. 2, we introduce the TDGL model and the $\omega$-gauge framework. We present the fully linearised mixed finite element scheme written in the $\omega$-gauge. In Sec. 3, we present our results in 2D. We first analyse a benchmark of the literature in a non convex domain and identify cases with numerical artefacts. Then, convergence orders are computing using the commonly used graphical method and the Richardson extrapolation technique. The analysis is continued with higher order finite elements. In Sec. 4, we extend our analysis to the 3D case and study three configurations: the unit cube, a sphere and a sphere with a geometrical defect.

 \section{The time-dependent Ginzburg-Landau model and the $\omega$-gauge framework}

\subsection{The time-dependent Ginzburg-Landau system}
The TDGL model describes the dynamics of a superconductor for temperatures close to the critical temperature \citep{gorkov1968,kato93} and is usually presented (in SI units) as

\begin{plain}
	\begin{equation}
		\eqalign{		
		\displaystyle	&\frac{\hbar^2}{2mD}\left( \frac{\partial}{\partial t} + i\frac{q}{\hbar}\phi\right) \psi= \frac{\hbar^2}{2m}\left( \nabla - i\frac{q}{\hbar}\mathbf{A} \right)^2\psi    - \alpha\psi - \beta|\psi|^2\psi,\cr
		\displaystyle &	\sigma \left( {\partial \mathbf{A} \over \partial t} + \nabla \phi \right) = {q \hbar \over 2 m i} (\psi^* \nabla \psi - \psi \nabla \psi^*) - {q^2 \over m}|\psi|^2\mathbf{A} - {1 \over \mu_0}\operatorname{\textbf{curl}} \left( \operatorname{\textbf{curl}} \mathbf{A} - \mu_0 \mathbf{H} \right),
	}
		\label{TDGLsi}
\end{equation}
\end{plain}
with  boundary conditions
\begin{plain}
	\begin{equation}
		\eqalign{		
		\displaystyle	\left( \nabla\psi - i\frac{q}{\hbar}\mathbf{A}\psi\right)\cdot \mathbf{n} = 0 &\mbox{ on } \partial \Omega, \cr
	\displaystyle	\left( \frac{1}{\mu_0}\operatorname{\textbf{curl}}\mathbf{A} \right)\times \mathbf{n} = \mathbf{H}\times \mathbf{n}  &\mbox{ on } \partial \Omega, \cr
	\displaystyle	\mathbf{E}\cdot\mathbf{n} = 0  &\mbox{ on } \partial \Omega, 
}
\label{TDGLbc}
\end{equation}
\end{plain}		
and initial conditions
\begin{plain}
	\begin{equation}
		\eqalign{		
		\displaystyle	\psi(\mathbf{x}, 0) = \psi_0(\mathbf{x}) &\mbox{ in }\Omega, \cr
			\displaystyle \mathbf{A}(\mathbf{x}, 0) = \mathbf{A}_0(\mathbf{x}) & \mbox{ in } \Omega. 
}
	\label{TDGLIc}
\end{equation}
\end{plain}

In previous equations, $\psi$ is the (complex valued) order parameter (with $\psi^*$ the complex conjugate) and $\mathbf{H}$  the applied magnetic field; $\alpha$ (negative) and $\beta$ (positive) are parameters depending on the temperature and the superconductor material; $q$ and $m$ denote the charge and the mass of the superconducting charge carrier, respectively. $D$ is a phenomenological diffusion coefficient and $\sigma$ has the dimension of an electrical conductivity. Finally, $\hbar$ is the reduced Planck constant and $\mu_0$ the magnetic permeability of the vacuum.

Numerical simulations are based on a non-dimensional form of \eqref{TDGLsi} which reads \citep{du1994}:

\begin{plain}
	\begin{equation}
		\eqalign{		
		\displaystyle&	\left( \frac{\partial}{\partial t} + i\kappa\phi\right) \psi=\left( \frac{1}{\kappa}\nabla - i\mathbf{A} \right)^2\psi + \psi -|\psi|^2\psi  \mbox{ in } \Omega,\cr
		\displaystyle& \left( \frac{\partial \mathbf{A}}{\partial t} +  \nabla\phi  \right) = - \operatorname{\textbf{curl}} \left(\operatorname{\textbf{curl}} \mathbf{A} - \mathbf{H}\right) +  \frac{1}{2i\kappa}(\psi^* \nabla \psi - \psi\nabla\psi^*) - |\psi|^2\mathbf{A}  \mbox{ in } \Omega,
}
	\label{TDGLScale2}
\end{equation}
\end{plain}

with boundary conditions
\begin{plain}
	\begin{equation}
		\eqalign{		
		\displaystyle	\left(	\frac{1}{\kappa} \nabla - i\mathbf{A} \right) \psi \cdot \mathbf{n} = 0 &\mbox{ on } \partial \Omega,\cr
		\displaystyle \operatorname{\textbf{curl}}\mathbf{A}\times \mathbf{n} = \mathbf{H}\times \mathbf{n} & \mbox{ on } \partial \Omega, \cr	\displaystyle	
		\mathbf{E}\cdot\mathbf{n} = 0 & \mbox{ on } \partial \Omega, 
}
	\label{TDGLbcScale}
\end{equation}
\end{plain}
and initial conditions
\begin{plain}
	\begin{equation}
		\eqalign{		
			\displaystyle \psi(\mathbf{x}, 0) = \psi_0(\mathbf{x}) & \mbox{ in }\Omega, \cr
			\displaystyle  \mathbf{A}(\mathbf{x}, 0) = \mathbf{A}_0(\mathbf{x}) & \mbox{ in } \Omega. 
		}
		\label{TDGLIcScaled}
	\end{equation}
\end{plain}

Lengths are scaled in units of the London penetration depth $\displaystyle\lambda = \left(\frac{m\beta}{\mu_0q^2(-\alpha)} \right)^\frac{1}{2}$. The ratio $\displaystyle \kappa = \frac{\lambda}{\xi}$, where $\displaystyle \xi = \frac{\hbar}{\sqrt{2m(-\alpha)}}$ is the coherence length, becomes the only physical parameter of the dimensionless formulation.

The non-dimensionalized Gibbs free energy ${\cal G}$ of the superconductor is \citep{du1992}:
\begin{equation}
	{\cal G}(\Psi, \mathbf{A}) = \int_\Omega \frac{1}{2}\left( |\psi|^2 - 1 \right)^2+ \left| \left( \frac{1}{\kappa} \nabla  - i \mathbf{A} \right) \psi \right|^2 + \left| \operatorname{\textbf{curl}} \mathbf{A} - \mathbf{H} \right|^2.
	\label{Gibbs}
\end{equation}

The time-dependent Ginzburg Landau equations \eqref{TDGLScale2} are related to the Gibbs energy through the following identities \citep{du2005}:

\begin{plain}
	\begin{equation}
		\eqalign{		
\displaystyle &\frac{\partial \psi}{\partial t} + i\kappa\phi\psi = -\frac{1}{2}\frac{\partial {\cal G}}{\partial \psi}\left( \psi, \mathbf{A}\right), \cr
 \displaystyle &\frac{\partial \mathbf{A}}{\partial t} + \nabla \phi = -\frac{1}{2}\frac{\partial {\cal G}}{\partial \mathbf{A}}\left( \psi, \mathbf{A}\right).
}
\label{tdglGibbRel}
\end{equation}
\end{plain}


\subsection{Gauge description}

Energy \eqref{Gibbs} is invariant under certain mathematical transformations called gauge transformations. Therefore, the physical properties of the system do not depend on these transformations.
In \cite{du1994} we find the general definition of a gauge for the TDGL model.

Given a function $\chi$, a gauge transformation is a linear transformation $G_\chi$ given by
\begin{plain}
	\begin{equation}
		\eqalign{	
			\displaystyle &	G_\chi(\psi, \textbf{A}, \phi) = (\zeta, \textbf{Q}, \Theta),\cr
			\displaystyle &	\mbox{where } \zeta = \psi \text{e}^{i\kappa\chi}, \quad\textbf{Q} = \textbf{A} + \nabla\chi, \quad \Theta = \phi -\frac{\partial \chi}{\partial t}. 
	}
		\label{GchiTDGL}
	\end{equation}
\end{plain}

Then $(\zeta, \textbf{Q}, \Theta)$ and $(\psi, \textbf{A}, \phi)$ solutions are said to be \emph{gauge equivalent}. It is easily seen from \eqref{GchiTDGL} that $\operatorname{\textbf{curl}}  \textbf{Q} = \operatorname{\textbf{curl}}  \textbf{A}$ and $|\zeta|^2 = |\psi|^2$. Hence the magnetic field or the density of the charge carriers, two physically relevant quantities, do not depend on the gauge.

For the definition of the $\omega$-gauge \citep{fleckinger97}, we define $\chi$ such that it satisfies the following boundary-value problem:
	\begin{align}		
			\displaystyle	&\left( \frac{\partial}{\partial t} - \omega\Delta \right) \chi  = \phi+ \omega\operatorname{div}\mathbf{A} \mbox{ in } \Omega \times (0, +\infty), \cr
			\displaystyle	&\omega\left( \mathbf{n} \cdot \nabla\chi \right)  =-\omega\left( \mathbf{n} \cdot \mathbf{A} \right) \mbox{ in } \partial \Omega \times (0, +\infty),
		\label{gaugeEq}
	\end{align}
with initial condition $\displaystyle \chi(\cdot, 0) = \chi_0$.
In this gauge, we have for $t > 0$:
\begin{align}\nonumber
		\displaystyle	\phi = -\omega\text{div}(\mathbf{A}),  &\mbox{ in } \Omega, \\
		\displaystyle	\omega\mathbf{A}\cdot \mathbf{n} = 0,  & \mbox{ in } \partial \Omega.
	\label{omegaG}
\end{align}
Each choice of $\omega$ corresponds to a different gauge: $\omega = 0$ gives the temporal  gauge, $\omega = 1$ the Lorenz gauge and $\omega =+\infty$ the Coulomb gauge. 

\subsection{The fully linearised MFE scheme}

In this section, we write the mixed variational formulation of the TDGL model under the $\omega$-gauge and the corresponding fully discretized scheme \citep{gaoSun2015, gaoSun2018}.

\subsubsection{Two-dimensional formulation}

In 2D, $\mathbf{A}$ has two components $A_1$ and $A_2$, depending on $x$ and $y$. As a result, the magnetic induction $\displaystyle \gamma = \operatorname{curl}\mathbf{A} = \frac{\partial A_2}{\partial x} - \frac{\partial A_1}{\partial y}$ is a scalar. Introducing $\gamma$ as a supplementary unknown, the system \eqref{TDGLScale2} with the $\phi = -\omega \text{div}(\mathbf{A})$ gauge can be rewritten as:
\begin{plain}
	\begin{equation}
		\eqalign{		 
			\displaystyle	&\frac{\partial \psi}{\partial t}  -i\kappa\omega\text{div}(\mathbf{A})\psi= 
			\left( \frac{1}{\kappa}\nabla -i\mathbf{A}\right)^2\psi+ \psi - |\psi|^2\psi,\cr
			\displaystyle &\gamma = \operatorname{curl} \mathbf{A},\cr
			\displaystyle	 &\frac{\partial {\mathbf A}}{\partial t} - \omega \nabla \text{div}(\mathbf{A}) + \operatorname{\textbf{curl}} \gamma
			=  \frac{1}{2i\kappa}\left( \psi^*\nabla\psi - \psi\nabla\psi^* \right)	- |\psi|^2{\mathbf A} + \operatorname{\textbf{curl}} {H},
		}
		\label{TDGLOmegaG}
	\end{equation}
\end{plain}
where $\displaystyle \operatorname{\textbf{curl}} = \left(\frac{\partial}{\partial x}, -\frac{\partial}{\partial y} \right)$. Boundary and initial conditions are:
	\begin{align}\nonumber
\frac{\partial \psi}{\partial{\mathbf n}} = 0, \quad \gamma= H, \quad \omega{\mathbf A}\cdot {\mathbf n} = 0  & \quad \mbox{ on } \partial\Omega \times (0, +\infty),\\
\psi({\mathbf x}, 0) = 1, \quad \gamma({\mathbf x}, 0)=0, \quad {\mathbf A}({\mathbf x}, 0) = (0,0)  & \quad \mbox{ on } \Omega.
\label{TDGLMFLorenzBC}
	\end{align}
To write the weak formulation, we introduce the following functional spaces
\begin{equation}
	\begin{array}{ll}
\displaystyle \text{H}^1 = \{ u \in L^2(\Omega), \nabla u \in \mathbf{L}^2(\Omega)\}, \\
	\displaystyle\mathbf{H}(\operatorname{div})=\left\{\mathbf{A} \mid \mathbf{A} \in \mathbf{L}^2(\Omega), \operatorname{div} \mathbf{A} \in L^2(\Omega)\right\},\\ \displaystyle\stackrel{\circ}{\mathbf{H}}(\text {div})=\left\{\mathbf{A} \mid \mathbf{A} \in \mathbf{H}(\text {div}),\left.\quad \mathbf{A} \cdot \mathbf{n}\right|_{\partial \Omega}=0\right\}.
	\end{array}
\label{defSobolevSpaces}
\end{equation}
We denote by  $\mathbf{H}$ (resp. $\cal H$), the Sobolev spaces corresponding to vector valued (resp. complex valued) functions. The dual space of a Sobolev space $\text{H}$ is denoted by $\text{H}^\prime$. The $L^2$ inner product is denoted by $(.,.)$.
The weak form corresponding to Eq. \eqref{TDGLOmegaG} is: find $\displaystyle \psi \in L^2(0,T;{\cal H}^1(\Omega))$ with $\displaystyle \frac{\partial\psi}{\partial t} \in L^2(0, T; {\cal H}^{-1}(\Omega))$ and $\displaystyle (\gamma, \mathbf{A}) \in L^2(0,T; \text{H}^1) \times \mathbf{L}^2(0,T; \overset{\circ}{\mathbf{H}}(\text{div}))$ with $\displaystyle \frac{\partial {\mathbf{A}}}{\partial t} \in \mathbf{L}^2(0, T; \overset{\circ}{\mathbf{H}}(\text{div})^\prime)$, where $\displaystyle \gamma = H$ on $\displaystyle \partial \Omega$, such that
\begin{plain}
	\begin{equation}
	\hspace{-5em}
	\eqalign{	
		\displaystyle & \left( \frac{\partial \psi}{\partial t}, w \right) - i\kappa {\omega} \left( (\operatorname{div}(\mathbf{A})\psi, w \right)
		= -\left( \left( \frac{1}{\kappa}\nabla -i\mathbf{A} \right)\psi, \left( \frac{1}{\kappa}\nabla -i\mathbf{A} \right) w \right) 
		+ \left( \left(  1 - |\psi|^2 \right) \psi, w \right) \quad \forall w \in {\cal H}^1(\Omega), \cr
		\displaystyle &\left(\gamma, \chi\right) - (\operatorname{\textbf{curl}} {\chi}, \mathbf{A}) = 0 \quad \forall \chi \in  \text{H}_0^1, \cr
		\displaystyle &\left( \frac{\partial \mathbf{A}}{\partial t}, \mathbf{v} \right)+ \left(\operatorname{\textbf{curl}} \gamma, \mathbf{v} \right)
		+ \left(  {\omega}\operatorname{div}\mathbf{A},  \operatorname{div}\mathbf{v}\right) - \frac{1}{2i\kappa} \left( \psi^*\nabla\psi - \psi\nabla\psi*, \mathbf{v}\right)  + (|\psi|^2\mathbf{A}, \mathbf{v}) = (\operatorname{\textbf{curl}} H, \mathbf{v})  \quad \forall \mathbf{v} \in \overset{\circ}{\mathbf{H}}(\text{div}),
	}
	\label{mixedVF}
\end{equation}
\end{plain}
for $t \in (0, T)$ with $\psi(x, 0) = \psi_0(x)$, $\mathbf{A}(x, 0) = \mathbf{A}_0(x)$ and $\gamma(x, 0) = \operatorname{curl}\mathbf{A}_0(x)$.\\ In numerical examples, we take $\mathbf{A}_0(\mathbf{x}) = (0,0)$ and $\psi_0(\mathbf{x}) = 1$ i.e. a pure superconducting state.

Following \cite{gaoSun2015}, we introduce the approximated fields $\mathbf{A}_h$, $\gamma_h$ and $\psi_h$ such that
\begin{equation}
	\displaystyle	{\mathbf A}_h \in {\mathbf U}_h^r, \quad
		\gamma_h \in V_h^{r+1}, \quad
		\psi_h \in V_h^r,
	\label{FEspacesMFETDGLLorenz}
\end{equation}
where $r \ge 0$. $U_h^r$ is the space of the Raviart-Thomas finite elements of order $r$ and $V_h^r$ the space of Lagrange finite elements of order $r$. In what follows, we omit the index $h$ for the fields.

The discrete formulation of \eqref{TDGLOmegaG} and \eqref{TDGLMFLorenzBC} is: find $\psi^{n+1}$ in $V_h^r$, $\gamma^{n+1}$ in ${V}_h^{r+1}$ and $\mathbf{A}^{n+1}$ in $U_h^r$, $r \geq 1$ such that for all $(w, \chi, \mathbf{v})$ in $
V_h^r \times \stackrel{\circ}{V}_h^{r+1} \times \stackrel{\circ}{{U}}_h^r$
\begin{plain}
	\begin{equation}
		\hspace{-8em}
		\eqalign{	
		\displaystyle	&\frac{1}{\delta t}(\psi^{n+1}, w) +\frac{1}{\kappa^2}(\nabla\psi^{n+1}, \nabla w) = \frac{1}{\delta t}(\psi^{n}, w) 
		+ \left( i\left(\kappa {\omega} + \frac{1}{\kappa}\right)\text{div}(\mathbf{A}^n)\psi^n, w\right) + \left( 2\frac{i}{\kappa}\psi^n\mathbf{A}^n, \nabla w \right) + ({\cal N}(\mathbf{A}^n, \psi^n), w), \cr
		\displaystyle	&(\gamma^{n+1}, \chi) - (\operatorname{\textbf{curl}} \chi, \mathbf{A}^{n+1}) = 0, \cr
		\displaystyle &	\frac{1}{\delta t}(\mathbf{A}^{n+1}, \mathbf{v}) + ({\omega}\text{div}(\mathbf{A}^{n+1}), \text{div}(\mathbf{v})) + (\operatorname{\textbf{curl}}  \gamma^{n+1}, \mathbf{v}) = 	\frac{1}{\delta t}(\mathbf{A}^{n}, \mathbf{v}) 	+ \frac{1}{2i\kappa}\left( \psi_n^*\nabla\psi_n - \psi_n\nabla\psi_n^*, \mathbf{v}\right) - ( |\psi^n|^2\mathbf{A}^n, \mathbf{v} )+ (\operatorname{\textbf{curl}} H, \mathbf{v}),	
}
	\label{rtDiscrete}
\end{equation}	
\end{plain}	
where ${\cal N}(\psi, \mathbf{A}) = \left(1 - \mathbf{A}^2 - |\psi |^2\right) \psi$.

\subsubsection{Three dimensional formulation}

In 3D, the TDGL model becomes:
\begin{plain}
	\begin{equation}
		\eqalign{	
			&\frac{\partial \psi}{\partial t}  -i\kappa\omega\text{div}(\mathbf{A})\psi= 
		\left( \frac{1}{\kappa}\nabla -i\mathbf{A}\right)^2\psi+ \psi - |\psi|^2\psi,\cr
		 & \boldsymbol{\gamma} = \operatorname{\textbf{curl}} \mathbf{A},\cr
			& \frac{\partial {\mathbf A}}{\partial t} - \omega \nabla \text{div}(\mathbf{A}) + \operatorname{\textbf{curl}} \boldsymbol{\gamma}
		=  \frac{1}{2i\kappa}\left( \psi^*\nabla\psi - \psi\nabla\psi^* \right)	- |\psi|^2{\mathbf A} + \operatorname{\textbf{curl}} {\mathbf H},
}
	\label{TDGL3DOmegaG}
\end{equation}
\end{plain}
with boundary and initial conditions

	\begin{align}\nonumber
	 	\frac{\partial \psi}{\partial{\mathbf n}} = 0, \quad \boldsymbol{\gamma} \times \mathbf{n}= \mathbf{H} \times \mathbf{n}, \quad \omega{\mathbf A}\cdot {\mathbf n} = 0 \quad & \mbox{ on } \partial\Omega \times (0, +\infty),\\
		\psi({\mathbf x}, 0) = 1, \quad \boldsymbol{\gamma}({\mathbf x}, 0)=0, \quad {\mathbf A}({\mathbf x}, 0) = (0,0) \quad & \mbox{ on } \Omega.
\label{TDGL3DMFLorenzBC}
\end{align}

To write the weak formulation, we introduce the following functional spaces:
\begin{equation}
	\begin{array}{ll}
	\displaystyle \mathbf{H}(\text {curl})=\left\{\mathbf{A} \mid \mathbf{A} \in \mathbf{L}^2(\Omega), \operatorname{\textbf{curl}} \mathbf{A} \in \mathbf{L}^2(\Omega)\right\},\\
	\displaystyle \stackrel{\circ}{\mathbf{H}}(\text {curl})=\left\{\mathbf{A} \mid \mathbf{A} \in \mathbf{H}(\text {curl}), \quad \mathbf{A} \times\left.\mathbf{n}\right|_{\partial \Omega}=\mathbf{0}\right\}.
\end{array}
\label{defSobolevSpaces2}
\end{equation}

The weak form corresponding to Eq. \eqref{TDGLOmegaG}  is: find $\displaystyle \psi \in L^2(0,T;{\cal H}^1(\Omega))$ with $\displaystyle \frac{\partial\psi}{\partial t} \in L^2(0, T; {\cal H}^{-1}(\Omega))$ and $\displaystyle (\boldsymbol{\gamma}, \mathbf{A}) \in \mathbf{L}^2(0,T; \mathbf{H}(\text {curl})) \times \mathbf{L}^2(0,T; \overset{\circ}{\mathbf{H}}(\text{div}))$ with $\displaystyle \frac{\partial {\mathbf{A}}}{\partial t} \in \mathbf{L}^2(0, T; \overset{\circ}{\mathbf{H}}(\text{div})')$, where $\displaystyle \boldsymbol{\gamma}\times\mathbf{n} = \mathbf{H}\times \mathbf{n}$ on $\displaystyle \partial \Omega$, such that
\begin{plain}
	\begin{equation}
		\hspace{-5em}
		\eqalign{	
		\displaystyle& \left( \frac{\partial \psi}{\partial t}, w \right) - i\kappa {\omega} \left( (\operatorname{div}\mathbf{A}\psi, w \right)
		=  -\left( \left( \frac{1}{\kappa}\nabla - i\mathbf{A} \right)\psi,\left( \frac{1}{\kappa}\nabla - i\mathbf{A} \right) w \right)+ \left( \left( 1 - |\psi|^2\right) \psi, w \right) \quad \forall w \in {\cal H}^1(\Omega), \cr
		\displaystyle &\left(\boldsymbol{\gamma}, \boldsymbol{\chi}\right) - (\textbf{curl }{\boldsymbol{\chi}}, \mathbf{A}) = 0 \quad \forall \boldsymbol{\chi} \in  \stackrel{\circ}{\mathbf{H}}(\text {curl}), \cr
		\displaystyle &\left( \frac{\partial \mathbf{A}}{\partial t}, \mathbf{v} \right)+ \left(\text{\textbf{curl} }\boldsymbol{\gamma}, \mathbf{v} \right)
		+ \left(  {\omega}\operatorname{div}\mathbf{A},  \operatorname{div}\mathbf{v}\right) - \frac{1}{2i\kappa} \left( \psi^*\nabla\psi - \psi\nabla\psi^*, \mathbf{v}\right) + (|\psi|^2\mathbf{A}, \mathbf{v}) = (\textbf{curl }\mathbf{H}, \mathbf{v})  \quad \forall \mathbf{v} \in \overset{\circ}{\mathbf{H}}(\text{div}),
}
	\label{mixedVF3D}
\end{equation}
\end{plain}
a.e. for $t \in (0, T)$ with $\psi(x, 0) = \psi_0(x)$,  $\mathbf{A}(x, 0) = \mathbf{A}_0(x)$ and $\boldsymbol{\gamma}(x, 0) = \textbf{curl }(\mathbf{A}_0(x))$. In  numerical examples, we take $\mathbf{A}_0(\mathbf{x}) = (0,0)$ and $\psi_0(\mathbf{x}) = 1$ i.e. a pure superconducting state.

To approximate the magnetic field $\boldsymbol{\gamma}$, we introduce the Nedelec FE space of order $r$ denoted by $Q_h^r$. The fully linearised scheme at the lowest order is: find $\psi^{n+1}$ in $V_h^1$, $\boldsymbol{\gamma}^{n+1}$ in $Q_h^0$ and $\mathbf{A}^{n+1}$ in $U_h^0$, such that for all $(w,\boldsymbol{\chi}, \mathbf{v})$ in $V_h^1 \times \stackrel{\circ}{\mathbf{Q}}_h^0 \times \stackrel{\circ}{\mathbf{U}}_h^0$
\begin{plain}
	\begin{equation}
		\hspace{-8em}
		\eqalign{	
		\displaystyle&	\frac{1}{\delta t}(\psi^{n+1}, w) +\frac{1}{\kappa^2}(\nabla\psi^{n+1}, \nabla w) = \frac{1}{\delta t}(\psi^{n}, w) 
		+ \left( i\left(\kappa {\omega} + \frac{1}{\kappa}\right)\text{div}(\mathbf{A}^n)\psi^n, w\right)  + \left( 2\frac{i}{\kappa}\psi^n\mathbf{A}^n, \nabla w \right) + \left({\cal N}(\mathbf{A}^n, \psi^n), w \right), \cr
		\displaystyle&	(\boldsymbol{\gamma}^{n+1}, \boldsymbol{\chi}) - (\textbf{curl } \boldsymbol{\chi}, \mathbf{A}^{n+1}) = 0 , \cr
		\displaystyle&	\frac{1}{\delta t}(\mathbf{A}^{n+1}, \mathbf{v}) + ({\omega}\text{div}(\mathbf{A}^{n+1}), \text{div}(\mathbf{v})) + (\textbf{curl } \boldsymbol{\gamma}^{n+1}, \mathbf{v}) = \frac{1}{\delta t}(\mathbf{A}^{n}, \mathbf{v})+ \frac{1}{2i\kappa}\left( \psi_n^*\nabla\psi_n - \psi_n\nabla\psi_n^*, \mathbf{v}\right) - ( |\psi^n|^2\mathbf{A}^n, \mathbf{v} ) + (\operatorname{\textbf{curl}} \mathbf{H}, \mathbf{v}).
}
	\label{rtDiscrete3D}
\end{equation}	
\end{plain}	

\pagebreak

\section{Convergence of the method in 2D}

\subsection{Benchmark in non convex geometry}

\begin{figure}[!h]
	\includegraphics[width = 0.25\textwidth]{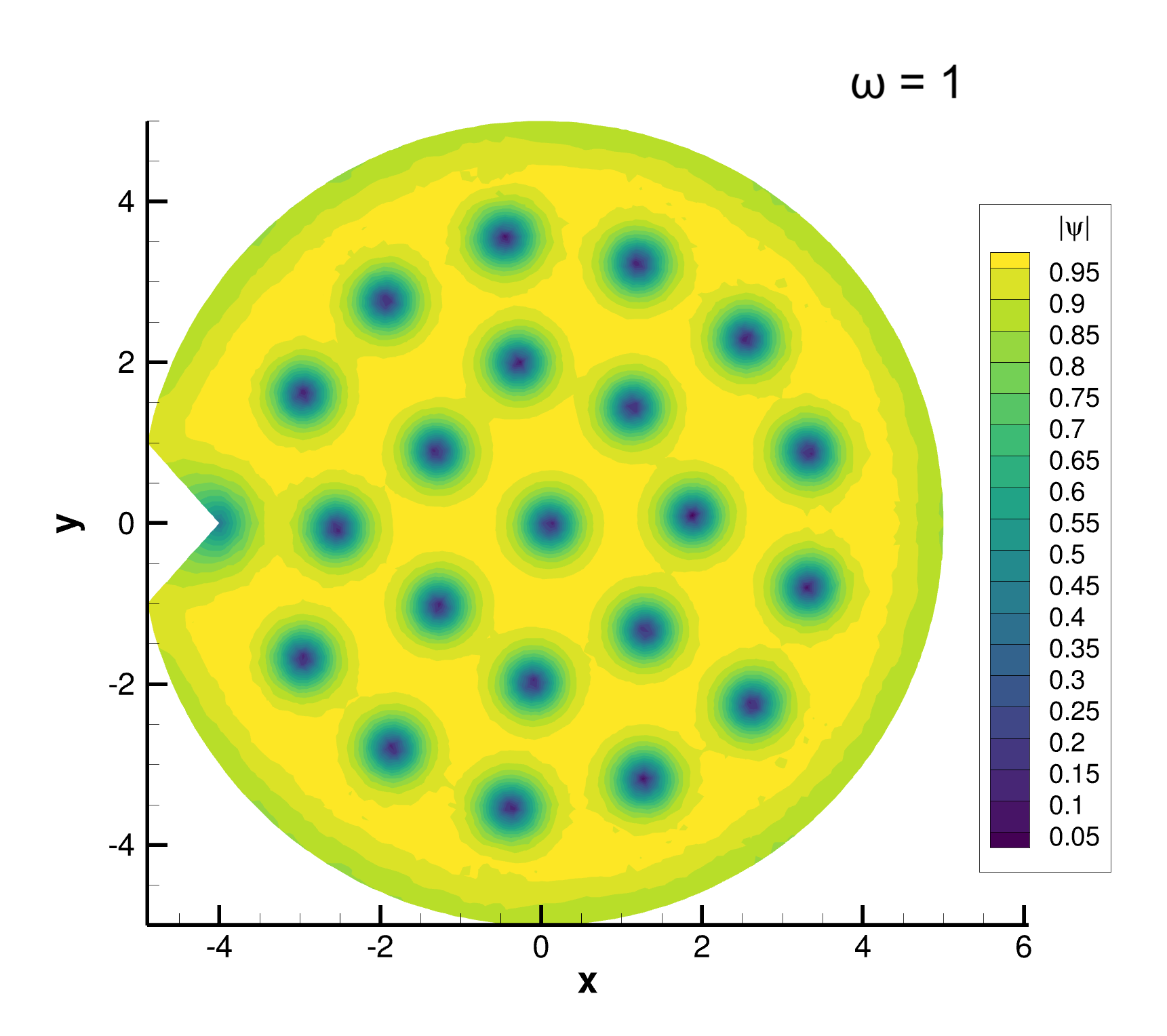}\hfill
	\includegraphics[width = 0.25\textwidth]{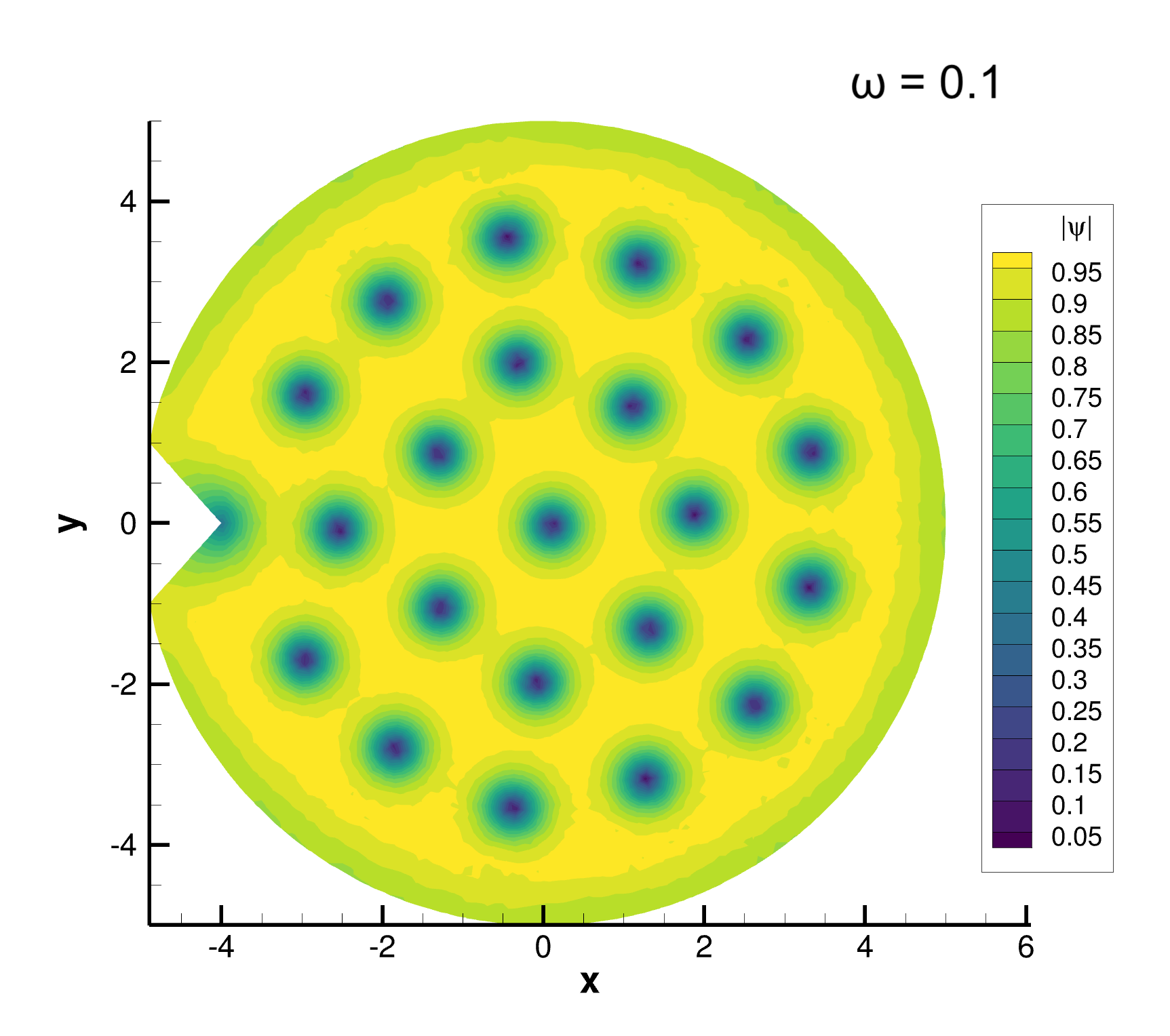}\hfill
	\includegraphics[width = 0.25\textwidth]{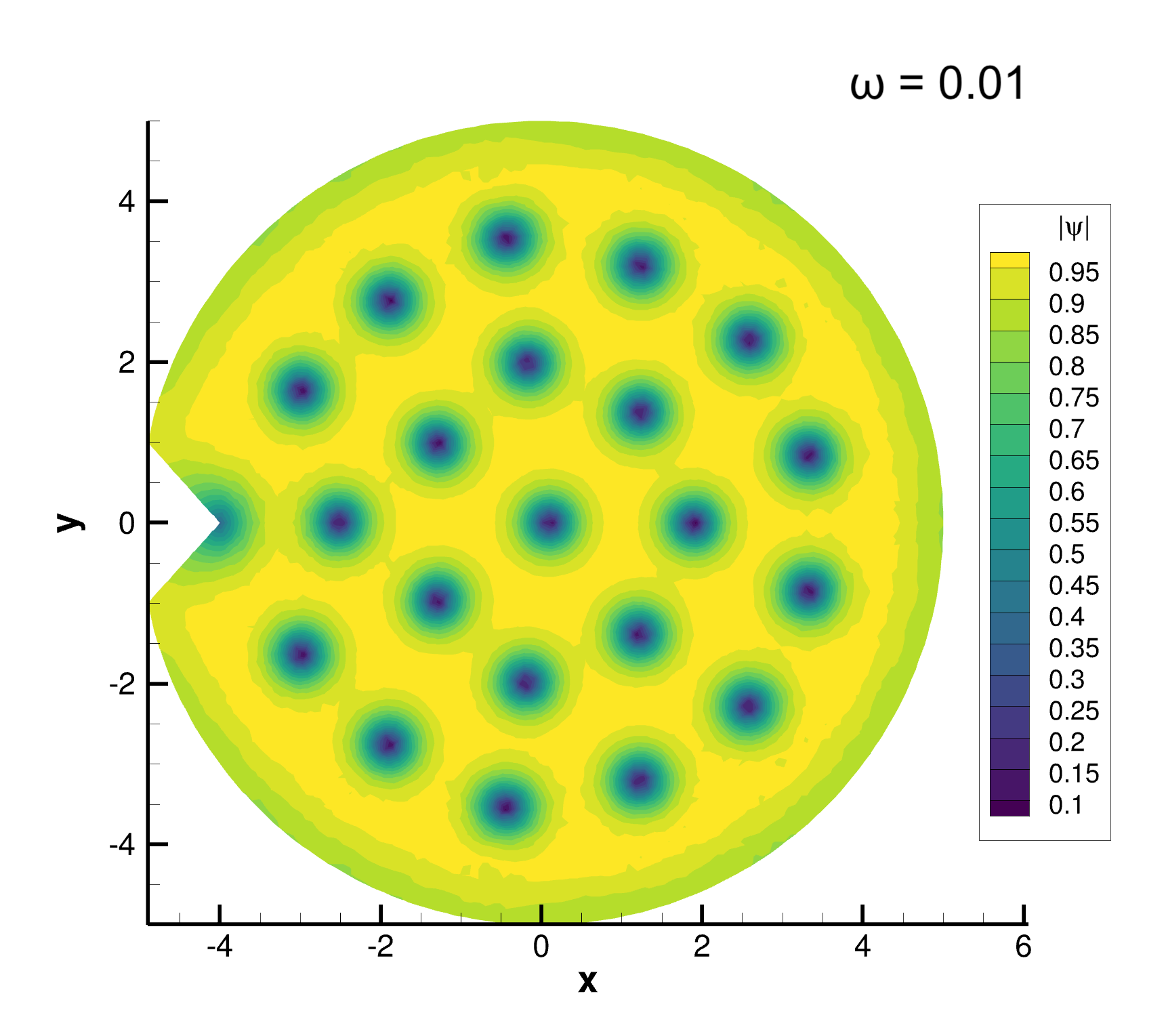}\\
	\includegraphics[width = 0.25\textwidth]{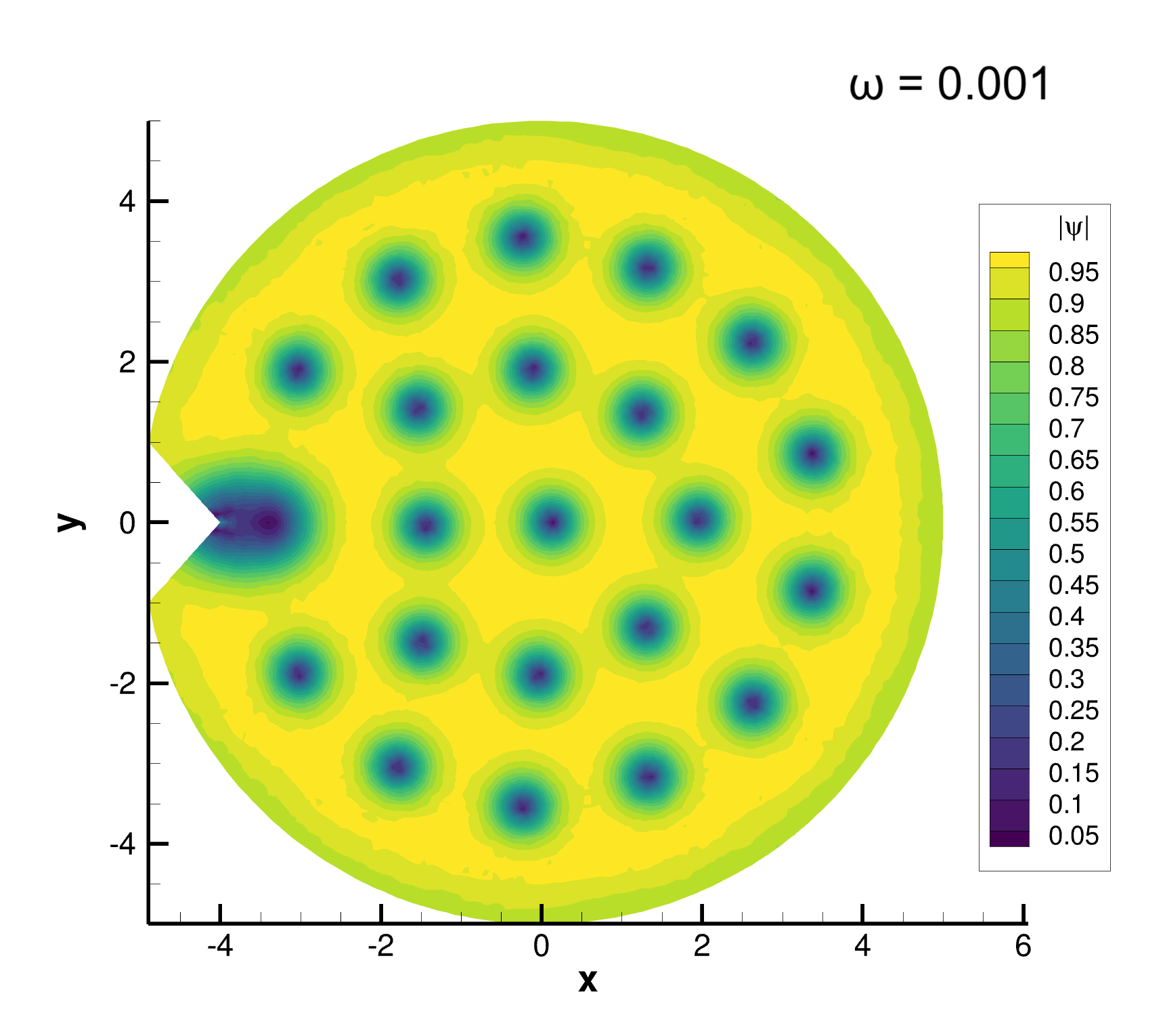}\hfill
	\includegraphics[width = 0.25\textwidth]{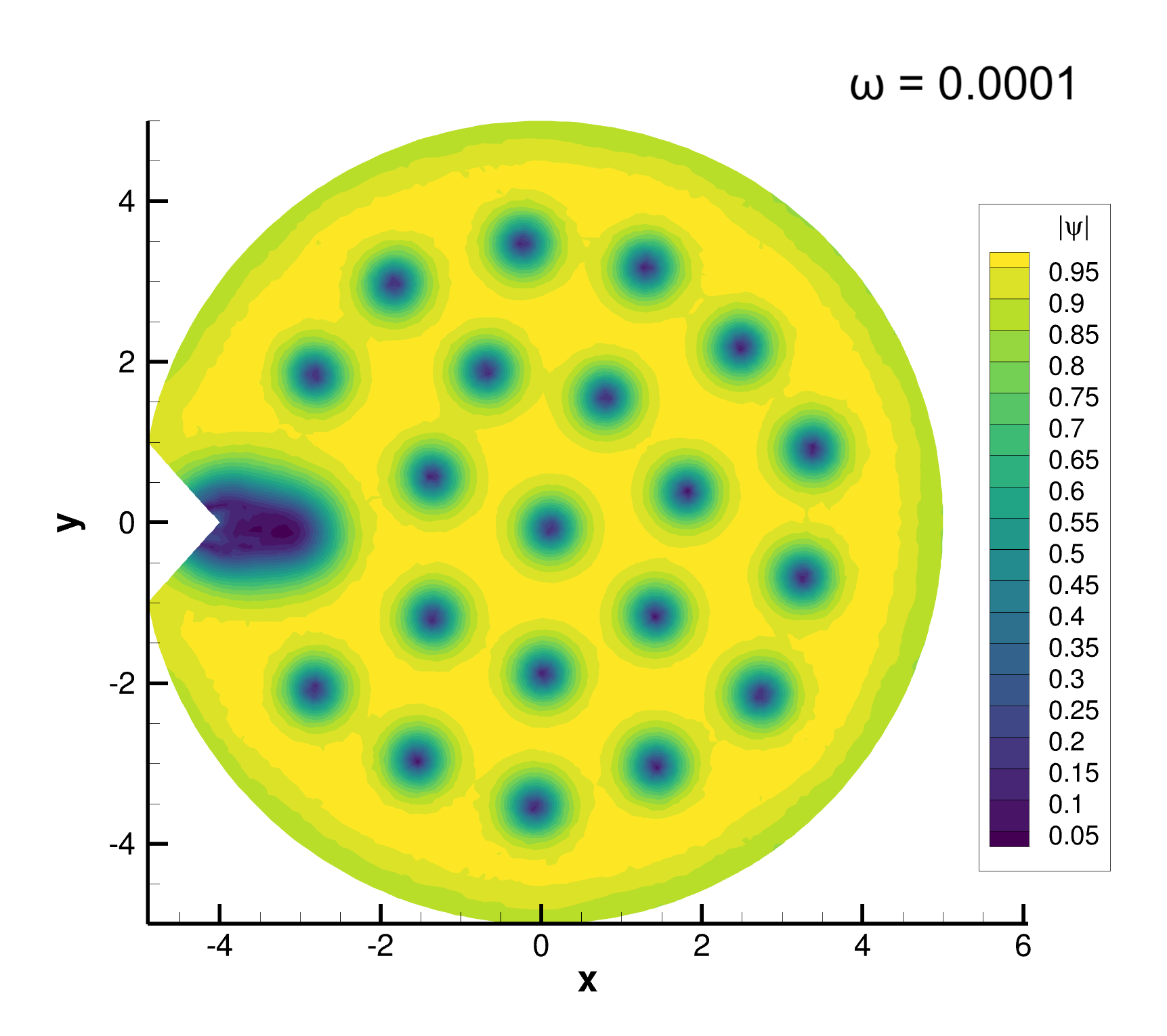}\hfill
	\includegraphics[width = 0.25\textwidth]{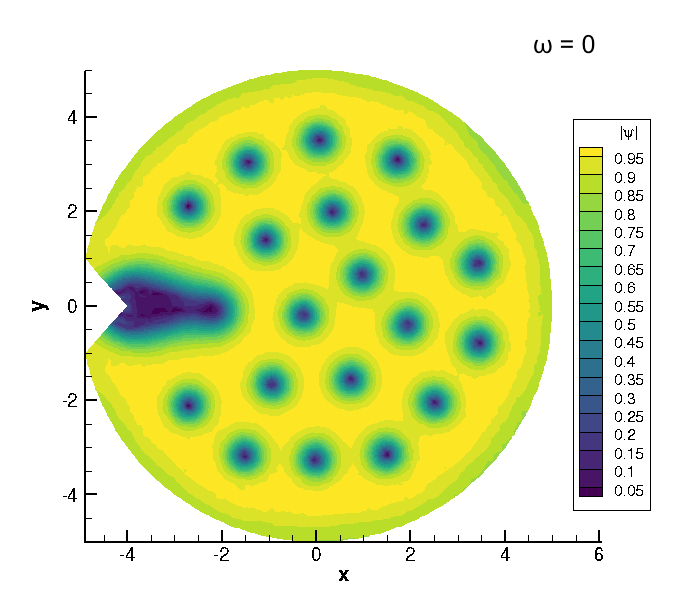}
	\caption{2D benchmark of a disk with an entrant corner. Finite elements of order r = 1 (see definition \eqref{FEspacesMFETDGLLorenz}). Contours of $|\psi|$ at $t = 5000$ for different values of the $\omega$ parameter.}
	\label{genGfig1}
\end{figure}

\begin{figure}[!h]
	\includegraphics[width = 0.25\textwidth]{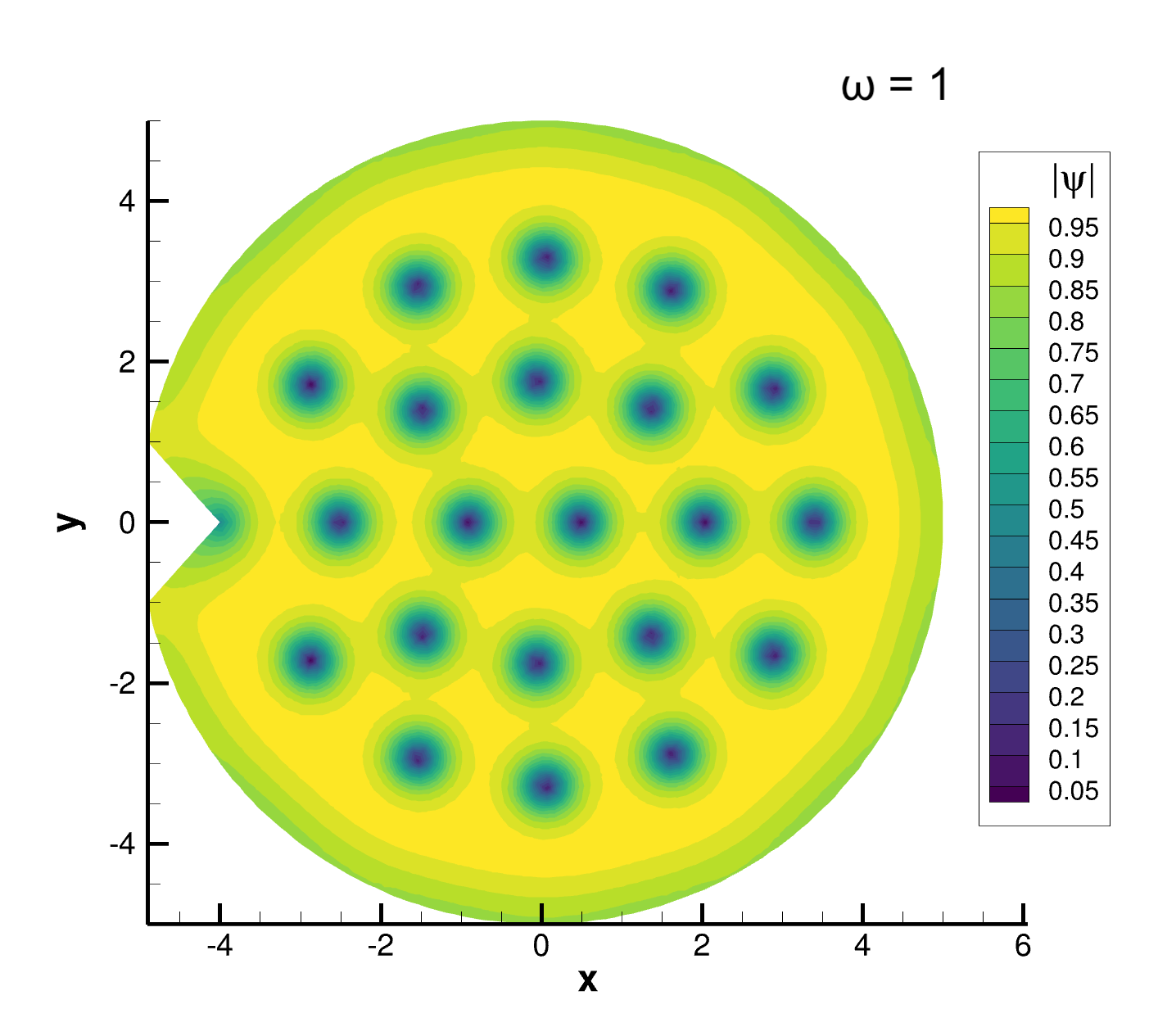}\hfill
	\includegraphics[width = 0.25\textwidth]{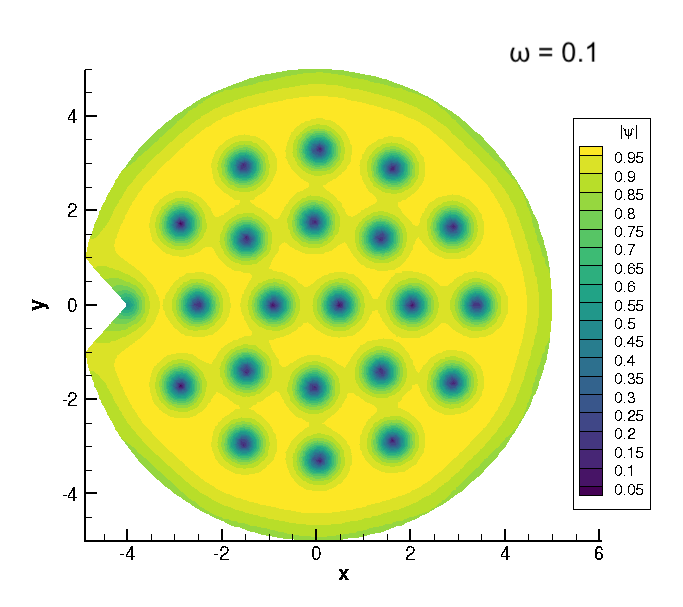}\hfill
	\includegraphics[width = 0.25\textwidth]{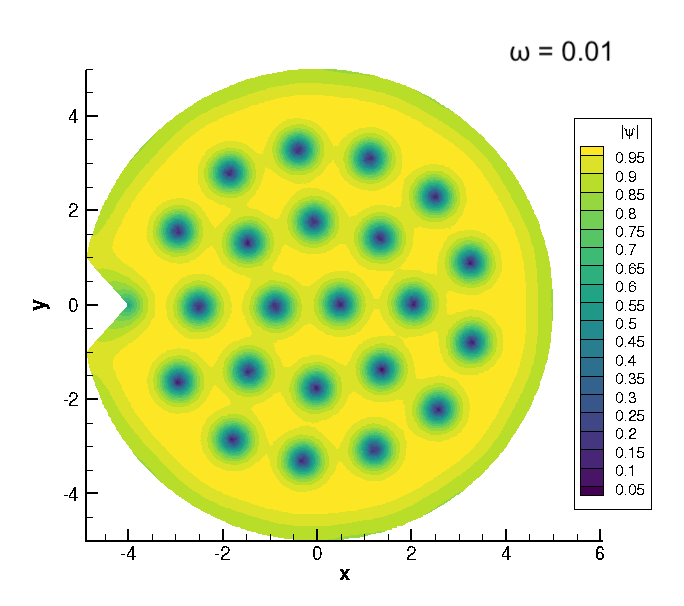}\\
	\includegraphics[width = 0.25\textwidth]{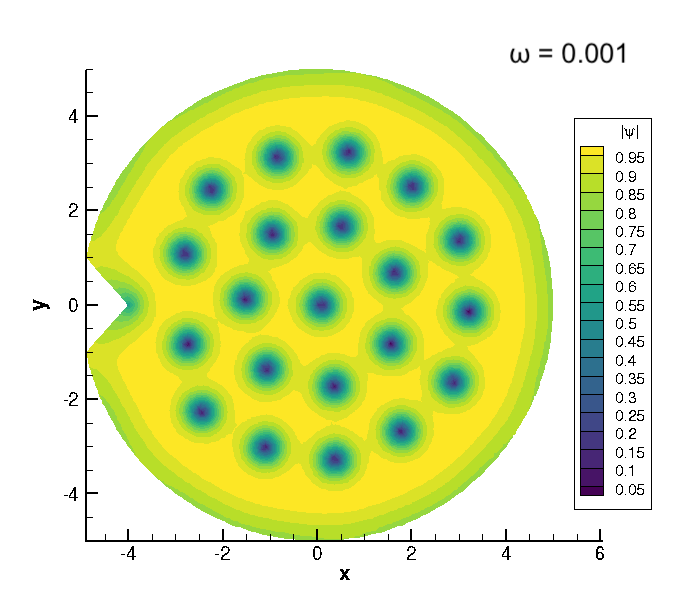}\hfill
	\includegraphics[width = 0.25\textwidth]{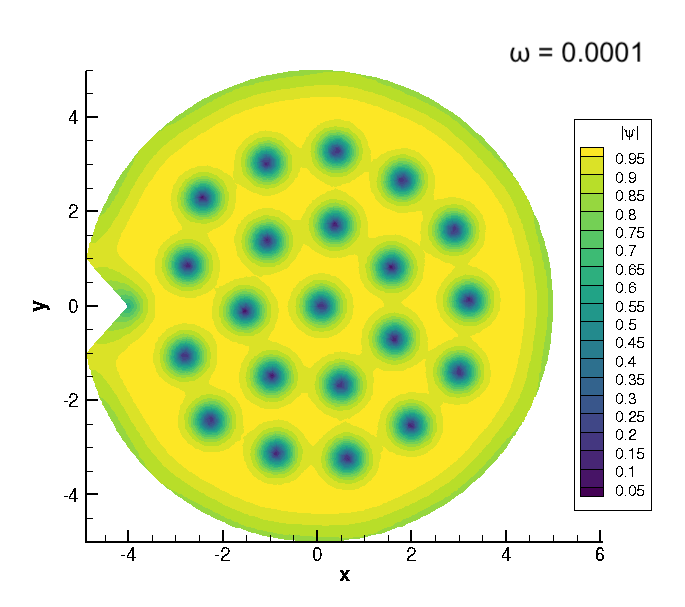}\hfill
	\includegraphics[width = 0.25\textwidth]{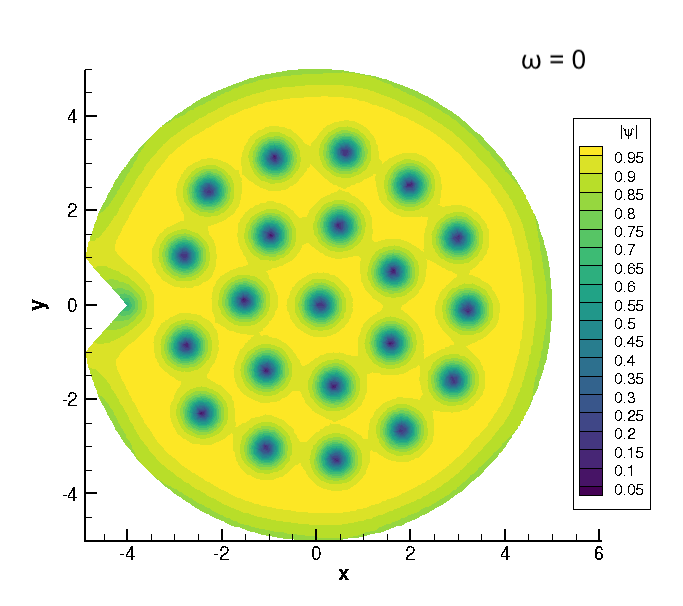}
	\caption{Same caption as Fig. \ref{genGfig1} with $r = 2$.}
	\label{genGfig2}
\end{figure}

In this section we study the TDGL model using the geometry of the disk with an entrant corner (see Fig. \ref{genGfig1}). This example was originally suggested in \cite{sorensen10}. We set $\kappa =4$, $H = 0.9$. The radius of the domain is $R = 5$ in units of $\lambda$. The mesh is uniform and the number of nodes per unit of the coherence length $\xi$ is 3. Using the mixed scheme \eqref{rtDiscrete}, Figs. \ref{genGfig1} and \ref{genGfig2} show the vortex pattern at $t =5000$ for FE of order $r = 1$ and $r =2$, respectively.

\newpage

For $r =1$ and $\omega \le 10^{-3}$, we notice a growing normal zone (i.e. a non-superconducting region where the order parameter takes very low values) near the indent. This zone might appear as an extended vortex \citep{sorensen10}, but in reality, it is just a numerical artefact. We can indeed resolve this zone by resorting to a finer mesh with 5 nodes per $\xi$. Normal zones are a common numerical issue with TDGL simulations \citep{Richardson04}.

For $r = 2$, there is no such normal region, but vortex arrangements at the final time are different. We observe 3 distinct vortex patterns. Table \ref{patternCaractericticDeg2} summarizes the characteristics of the final state for  each $\omega$ when $r = 2$. ${\cal G}_n$ is the free energy computed at time $t = n\delta t$.

\begin{table}[h!]
	\centering
	\begin{tabular}{lclclclclclcl}
		\hline
		$\omega$& Number of vortices &$|{\cal G}_{n} - {\cal G}_{n-1} |$ &${\cal G}_{n_{max}}$ &$n_{max}$\\
		& & at $n = n_{max}$&&&\\
		\hline
		$1$&21& $< 10^{-10}$&16.4711& 5000\\
		\hline
		$10^{-1}$&21& $< 10^{-10}$&16.4711&5000\\
		\hline
		$10^{-2}$&22&$1.6\cdot 10^{-7}$ &16.0959&5000\\
		\hline
		$10^{-3}$&21&$7.4\cdot 10^{-8}$ &16.4362&5000\\
		\hline
		$10^{-4}$&21&$4.3\cdot 10^{-7}$ &16.4310&5000\\
		\hline
		$0$&21&$1.3\cdot 10^{-6}$ &16.4338&5000\\
		\hline
	\end{tabular}
	\caption{2D benchmark of a disk with an entrant corner. Finite elements of order r = 2. Characteristics of vortex patterns for each $\omega$ at $t = n_{max}\delta t$ with $n_{max} = 5000$.}
	\label{patternCaractericticDeg2}
\end{table}

\pagebreak

\begin{figure}[!h]
	\centering
	\includegraphics[width = 0.5\textwidth]{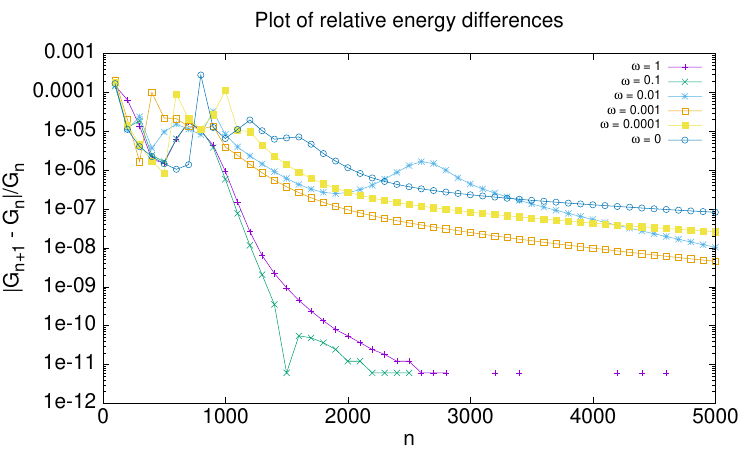}
	\caption{2D benchmark of a disk with an entrant corner. Finite elements of order r = 2. Relative energy difference $|{\cal G}_{n+1} -{\cal G}_n|/{\cal G}_n$ in logarithmic scale for $\omega = 1,10^{-1}, 10^{-2}, 10^{-3}, 10^{-4}, 0$;}
	\label{genGfig4}	
\end{figure}

Figure \ref{genGfig4} shows the relative energy differences $\displaystyle \frac{|{\cal G}_{n+1} - {\cal G}_n|}{{\cal G}_n}$ for $n = 0 \dots 5000$ for $w = 1, 10^{-1}, 10^{-2}, 10^{-3}, 10^{-4}, 0$ for the case $r = 2$. We observe that cases $\omega = 1, 10^{-1}$ are the fastest cases for reaching the equilibrium. We also notice that the case $\omega = 10^{-2}$ converges faster than cases with lower values of $\omega$. 
 
\begin{figure}[!h]
	\includegraphics[width = 0.4\textwidth]{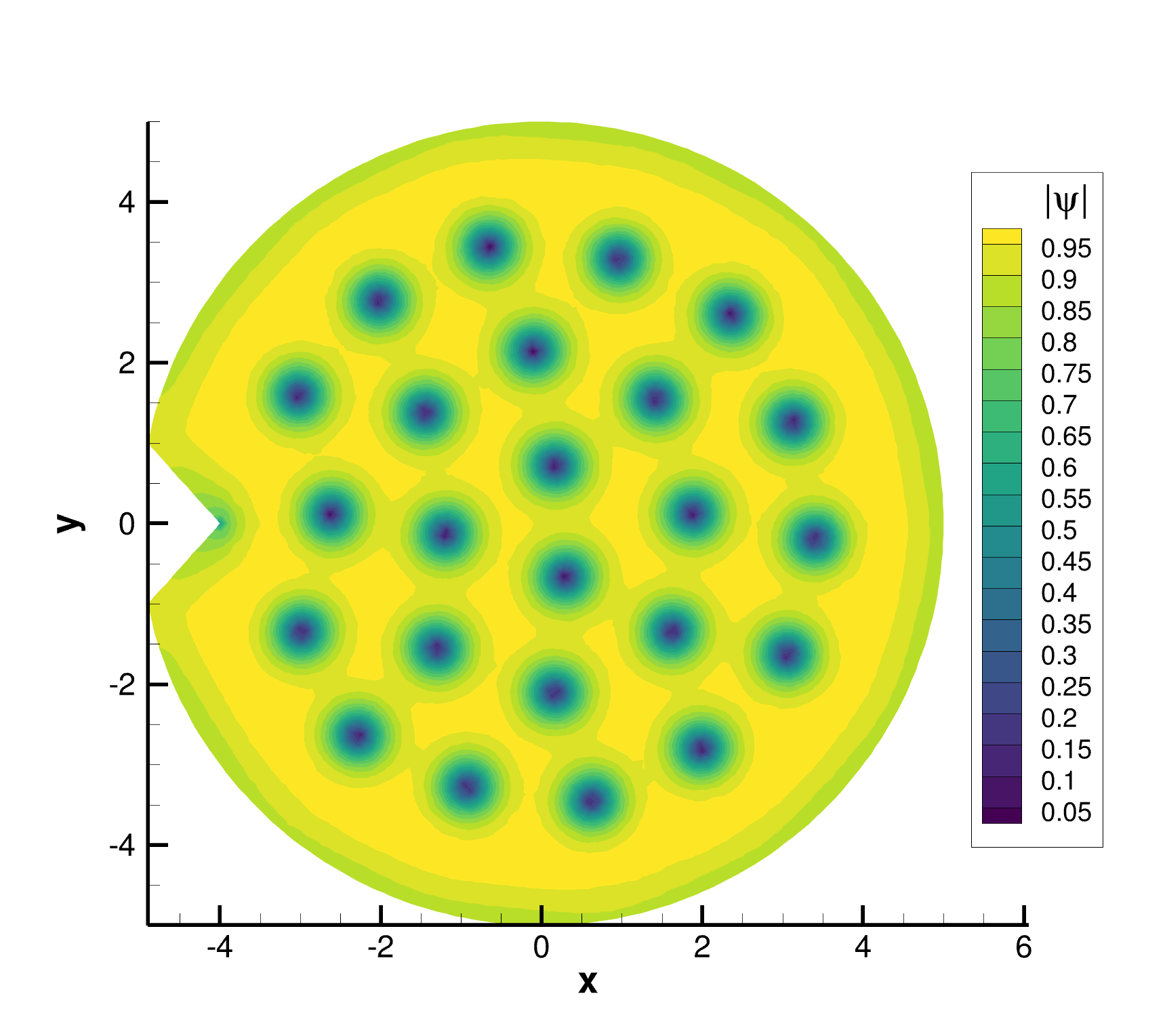}\hfill
	\includegraphics[width = 0.4\textwidth]{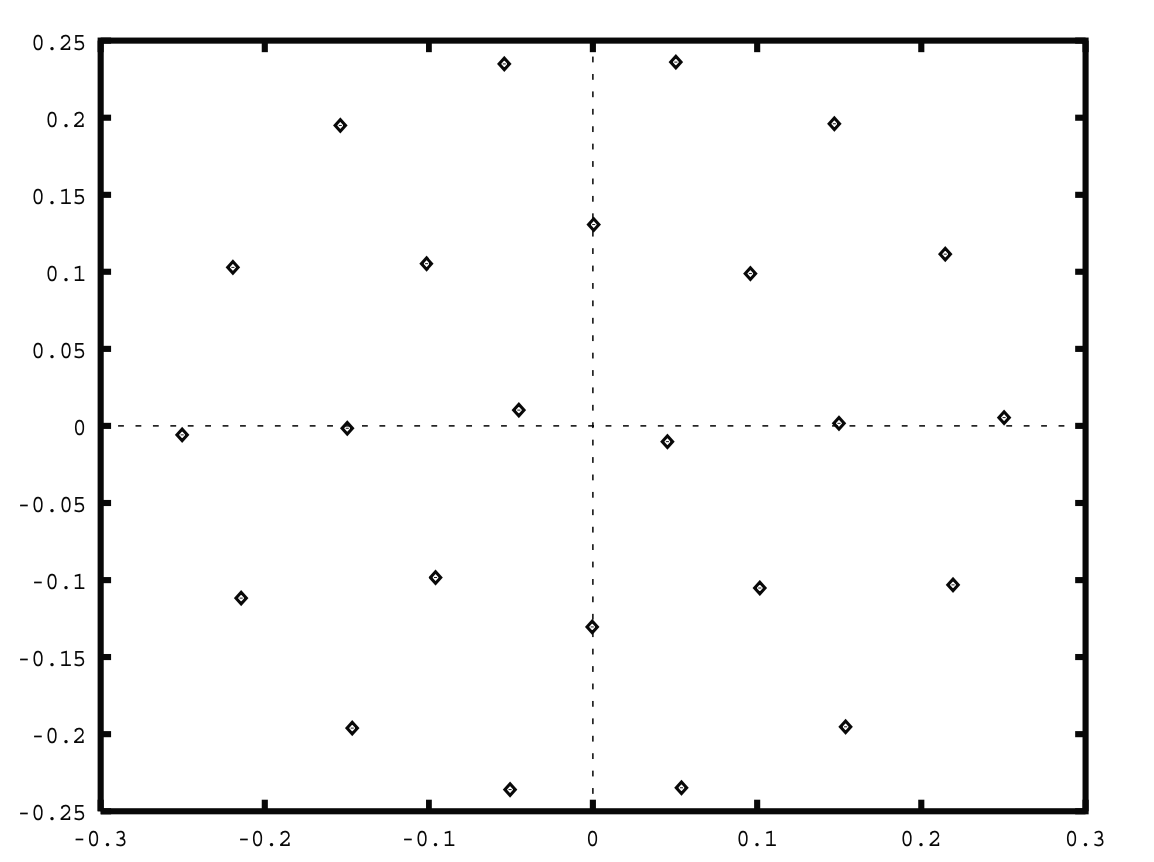}
	\caption{2D benchmark of a disk with an entrant corner. Finite elements of order r = 2. Configuration of lowest energy with 24 vortices (left). Minimizer of \eqref{serfatyWn} corresponding to a system of 24 point vortices, taken from \cite{gueron1999discrete}(right).}
	\label{configLowEnerg}
\end{figure}

Note that vortex arrangements could be different, depending on the value of $\omega$. Each state corresponds to a numerically found local minimizer of the energy \eqref{Gibbs}. Among these minimizers, the ground state is defined as the global minimum. Table \ref{patternCaractericticDeg2} shows that the final state corresponding to $\omega = 10^{-2}$ has the lowest energy.
In Fig. \ref{configLowEnerg} (left panel), we show another configuration with 24 vortices and energy equal to $15.5547$. It has been obtained by starting with the final state corresponding to $\omega = 10^{-4}$ and then progressively increasing the gauge parameter up to $\omega = 1$. The vortex pattern does not have a symmetry with respect to the $x$-axis unlike the ones found hitherto.  This configuration corresponds well to that numerically found in 
 \cite{gueron1999discrete} (see the right panel in Fig. \ref{configLowEnerg}) as a minimizer with also $n=24$ vortices, but for the renormalized energy of a system of $n$ point vortices  \cite{sandier2008vortices}:
\begin{equation}
	\displaystyle	w_n(x_1, \dots, x_n) =-\pi\underset{i\ne j}{\sum}\text{log}|x_i-x_j| + C\pi n \sum_{i=1}^{n}|x_i|^2.
	\label{serfatyWn}
\end{equation}

To conclude this part, we note that the convergence towards the equilibrium is faster when $\omega \geq 10^{-2}$.
In the following sections, we compute convergence orders for different choices of the gauge.

\subsection{A manufactured TDGL system}

A manufactured system is a system for which the exact solution is known analytically.  The general idea 
of the technique of manufactured solutions (e. g. \cite{BEC-book-1998-Roache}) is to modify the
original system of equations by introducing an extra source term, such that the new system admits an exact solution
given by a convenient analytic expression. Even though in most cases exact solutions constructed in this way are not physically realistic,
this approach allows one to rigorously verify computations.

In the case of the TDGL system, the manufactured system on the unit square $(0,1)\times(0,1)$ is
\begin{plain}
	\begin{equation}
		\eqalign{
		&	\frac{\partial \psi}{\partial t}  -i\kappa\omega\text{div}(\mathbf{A}) 
		-\left( \frac{1}{\kappa}\nabla -i\mathbf{A}\right)^2\psi - \psi +|\psi|^2\psi = g,\cr
		&	 \frac{\partial {\mathbf A}}{\partial t} - \omega \nabla \text{div}(\mathbf{A}) + \operatorname{\textbf{curl}}  \operatorname{{curl}} \mathbf{A}
		- \frac{1}{2i\kappa}\left( \psi^*\nabla\psi - \psi\nabla\psi^* \right)	+ |\psi|^2{\mathbf A} =\operatorname{\textbf{curl}} {H}+ \mathbf{f} ,
}
	\label{TDGLGeneralizedGArtificial}
\end{equation}
\end{plain}

with boundary and initial conditions
\begin{equation}
	\begin{array}{ll}
		\nabla\psi \cdot \mathbf{n} = 0, \quad \operatorname{{curl}} \mathbf{A} = H, \quad \omega\mathbf{A}\cdot \mathbf{n} = 0, \quad \mbox{on } \partial\Omega,
	\end{array}
	\label{TDGLBcGeneralizedGArtificial1}
\end{equation}	
where $\mathbf{f}$ and $g$ are defined such that the exact solution of \eqref{TDGLGeneralizedGArtificial} reads:
\begin{plain}
	\begin{equation}
		\eqalign{
	&	\psi = \exp(-t)\left( \cos(\pi x) + i\cos(\pi y)\right), \cr
	&	\mathbf{A} = 
	\left(
\begin{array}{ll}
			\exp(y-t)\sin(\pi x) \\
			\exp(x-t)\sin(\pi y)
		\end{array}\right),\cr
	&	H = \exp(x-t)\sin(\pi y) - \exp(y - t) \sin(\pi x).
}
	\label{artificalSol}
\end{equation}
\end{plain}

It is shown in \cite{gaoSun2018} that, if  the exact solution of \eqref{mixedVF3D} is regular enough, then for a final time $t_N$ of the scheme, the following error estimates hold:
\begin{plain}
	\begin{equation}
		\eqalign{
		\displaystyle	||\psi_h^N -\psi^N||_{L_2} = O(\Delta t + \Delta x^{r+1}),\cr
		\displaystyle	||\mathbf{A}_h^N -\mathbf{A}^N||_{L_2} = O(\Delta t + \Delta x^{r+1}),\cr
		\displaystyle	\Delta t \sum_{n=1}^{N} ||\gamma_h^n - \operatorname{\textbf{curl}} \mathbf{A}^n||^2_{L_2} = O(\Delta t^2 + \Delta x^{2r+2}).
}
	\label{errGao2015}
\end{equation}
\end{plain}

\subsection{Convergence analysis of the scheme \eqref{rtDiscrete} for the case $r = 1$}

We now describe two methods to compute the convergence orders. The graphical method is the usual technique. However it becomes computationally costly for higher order finite elements. Therefore, we use the Richardson extrapolation method that was proved to be fast and reliable.

\subsubsection{The graphical method}
\begin{figure}[!h]
	\centering
	\includegraphics[width = 0.48\textwidth, angle = 0]{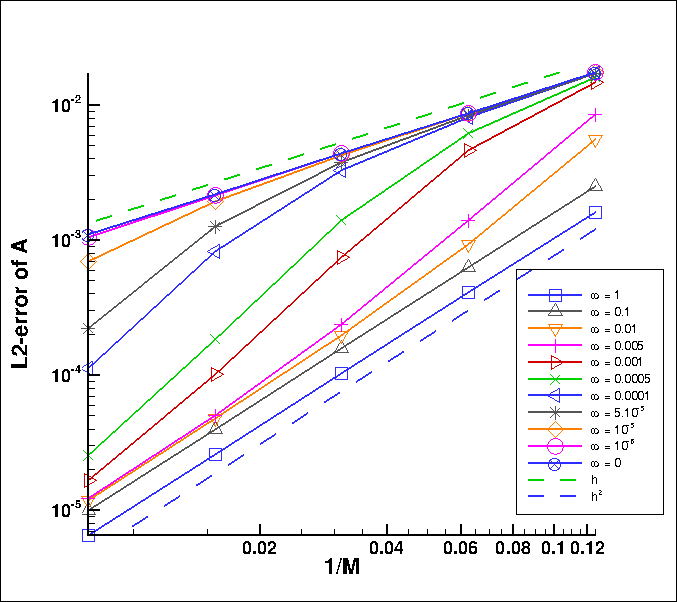}\hfill
	\includegraphics[width = 0.48\textwidth, angle = 0]{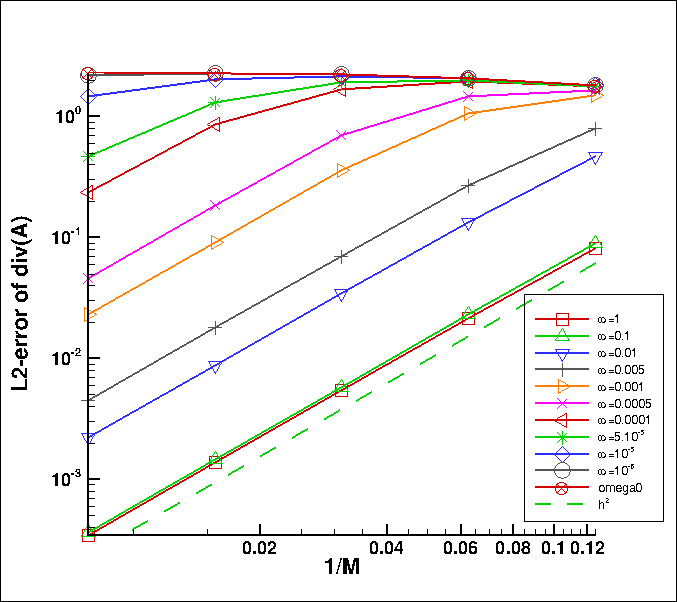}
	\caption{2D manufactured solutions benchmark. Finite elements of order r = 1. Orders for the vector potential $\mathbf{A}$ (left) and its divergence $\text{div}(\mathbf{A})$ (right).}
	\label{orderGeneralizedG}
\end{figure}

\begin{figure}[!h]
	\centering
	\includegraphics[width = 0.3\textwidth, angle = 0]{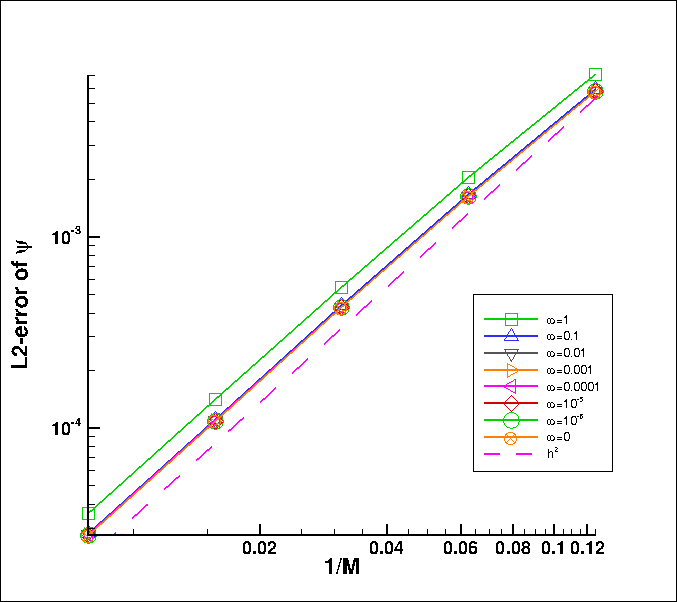}\hfill
	\includegraphics[width = 0.3\textwidth, angle = 0]{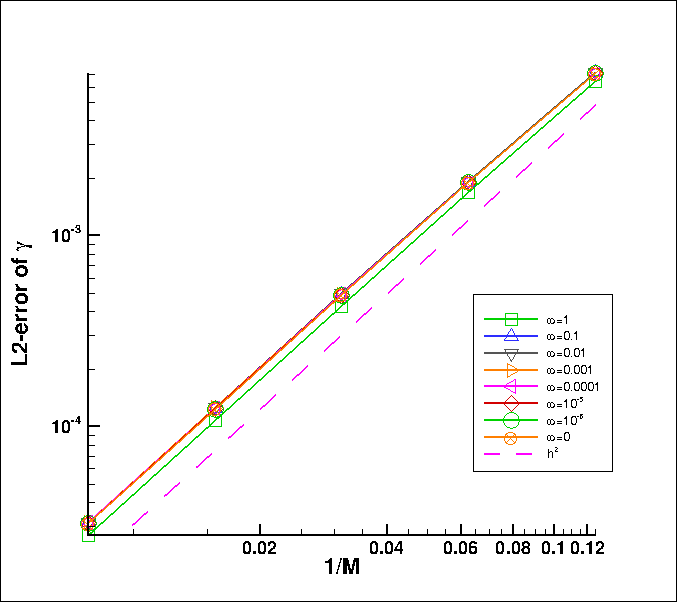}\hfill
	\includegraphics[width = 0.3\textwidth, angle = 0]{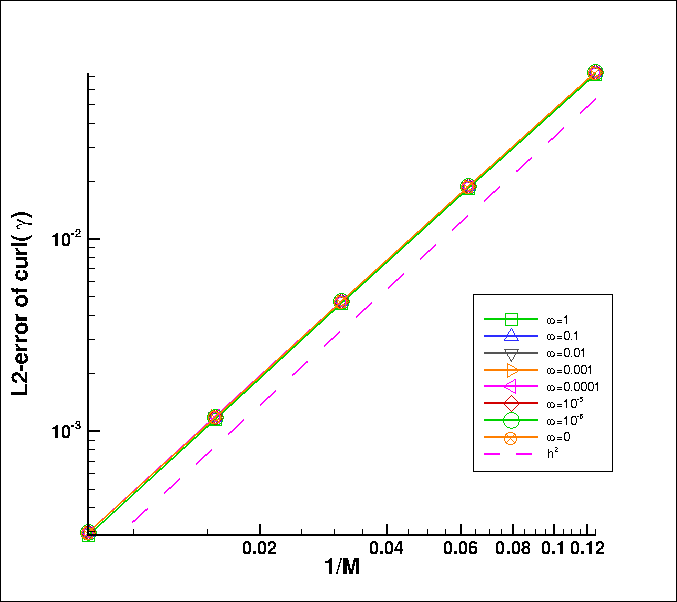}
	\caption{2D manufactured solutions benchmark.  Finite elements of order r = 1. Orders for  $\psi$, $\gamma= \operatorname{\textbf{curl}} \mathbf{A}$ and $\operatorname{\textbf{curl}}\gamma$. }
	\label{orderGeneralizedG1}
\end{figure}

We choose $\displaystyle \Delta x = \frac{1}{M}$ and $\displaystyle \Delta t = \frac{1}{M^{3}}$ with $M = 8,16,32,64,128$ and iterate the scheme $\displaystyle \frac{M^3}{8}$ times. We compute the solution at $\displaystyle t = \frac{1}{8} = 0.125$ and then compare it with the exact solution \eqref{artificalSol}. Results are shown in Figs.  \ref{orderGeneralizedG} - \ref{orderGeneralizedG1} with $M = 8,16,32,64,128$.

We observe that the vector potential loses one order between $\omega = 10^{-2}$ and $\omega = 10^{-5}$.  Between $\omega = 10^{-3}$ and $\omega = 5\times10^{-5}$ we observe an increase of the order for $\mathbf{A}$ at smaller sizes of the mesh; this suggests the existence of an inflection point where the order is maximum;  this point depends on the size of the mesh. As regards the divergence of $\mathbf{A}$, it loses two orders between $\omega = 10^{-2}$ and $\omega = 10^{-5}$; between $\omega = 10^{-3}$ and $\omega = 10^{-5}$, the convergence rate decreases monotonously as the mesh size increases; the beginning of the decrease depends on the size of the mesh.
The orders for $|\psi|$, $\gamma$ and $\operatorname{\textbf{curl}}\gamma$ are not affected. 

The graphical method used to determine the orders can be misleading, since the relation between $\Delta t$ and $\Delta x$ is fixed and imposed by the expected convergence rate.  Indeed, we are not able to see orders greater than $q$ if $\displaystyle \Delta t = \frac{1}{M^q}$. Besides the method is time consuming, since we need to compute a number of iterations of order $N = \text{O}(M^q)$.

\subsubsection{The Richardson extrapolation method}
 In this section, we describe another method to estimate convergence rates, known as the Richardson extrapolation technique. The advantage of this method is that the number of iterations $N_{iter}$ is fixed. It is very efficient, once the parameters $\Delta t$ and $N_{iter}$ are chosen carefully. 

Let us first describe the method. Consider a quantity $u$ to be evaluated numerically. We denote by $u_{ex}$ its exact value and $u_h$ its approximation, where $h$ is the discretization step to be refined. If $p$ is the order of the numerical scheme, then we assume that:
\begin{equation}
	\begin{array}{ll}
		\displaystyle	u_{ex} = u_h + Ch^p, \\
		\displaystyle	u_{ex} = u_{\frac{h}{2}} + C\left( \frac{h}{2} \right)^p, \\
		\displaystyle u_{ex} = u_{\frac{h}{4}} + C\left( \frac{h}{4} \right)^p.	
	\end{array}
	\label{richardson}
\end{equation}
By substitution we deduce that the convergence order can be computed as:
\begin{equation}
	\displaystyle	p =\frac{1}{\log 2} \log \left( \frac{u_{\frac{h}{2}} - u_h}{u_{\frac{h}{4}} -u_{\frac{h}{2}}} \right).
	\label{richardonOrder}
\end{equation}
Since $u$ is a field in our case, we consider $\displaystyle p =\frac{1}{\log 2} \log \left( \frac{ ||u_{\frac{h}{2}} - u_h||_{L^2}}{ ||u_{\frac{h}{4}} -u_{\frac{h}{2}}||_{L^2}} \right)$.

As a benchmark, we first recover the results in the case of the Lorenz gauge.
Table \ref{table:1} shows the orders obtained in this case for $M = 16$, so that $\displaystyle \Delta x =\frac{1}{16}= 0.0625$. The time step is $\Delta t = 10^{-3}$ and we take $N_{iter} = 125$ iterations. We recover the correct orders.

\begin{table}[h!]
	\begin{tabular}{lclclclclclcl}
		\hline
		&  $\text{Err}_\psi$ &$\text{Err}_\mathbf{A}$ &$\text{Err}_\gamma$& $\text{Err}_{\operatorname{\textbf{curl}}\gamma}$&$\text{Err}_{\text{div}(\mathbf{A})}$ \\
		\hline
		$||u_{\frac{\Delta x}{2}} - u_{\Delta x}||_{L^2} $& $0.0029113$ & $0.00145666 $& $0.00129579 $ & $0.00987884$&$ 0.00686983$\\
		\hline
		$||u_{\frac{\Delta x}{4}} -u_{\frac{\Delta x}{2}}||_{L^2}$ & $0.000729874$ & 	$0.00036462$9 & $0.00032392$ & $0.00245809$& $0.001727$\\
		\hline
		$\text{order} $&  $1.99594$ & $1.9982$ & $2.00008$  &  $2.00681$ &$  1.99169$\\
		\hline
	\end{tabular}
	\caption{2D manufactured solutions benchmark. Finite elements of order r = 1. Computed $L^2$ norm errors of $\psi$, $\mathbf{A}$, $\gamma$, $\operatorname{\textbf{curl}}\gamma$ and $\text{div}(\mathbf{A})$ for the Lorenz gauge.}
	\label{table:1}
\end{table}

Table \ref{tableDeg1allOm} shows the orders computed for 10 values of $\omega$ between $0$ and $1$ with  $M = 16$, $\Delta t = 10^{-3}$ and $N_{iter} = 125$. The orders for quantities $|\psi|$, $\gamma$ and $\operatorname{\textbf{curl}}\gamma$ are in agreement with the graphical results of Fig. \ref{orderGeneralizedG1}. In Tabs. \ref{tableDeg1VectPotCompGraphVsRich} - \ref{tableDeg1DivACompGraphVsRich}, we compare the two methods for $\mathbf{A}$ and $\operatorname{div}\mathbf{A}$. The reported graphical values are the average slopes obtained with the 3 points used in the Richardson method; in this case they correspond to $M = 16, 32,64$. We observe a good agreement between the two techniques.

\begin{table}[h!]
	\hspace{-6em}
	\footnotesize
	{
		\begin{tabular}{lclclclclclcl}
			\hline
			&  $\omega = 1$ & $\omega = 10^{-1}$ & $\omega = 10^{-2}$ & $\omega = 10^{-3}$ & $\omega = 10^{-4}$ & $\omega = 10^{-5}$ &$\omega = 10^{-6}$ & $\omega = 0$\\ 
			\hline
			$\text{Err}_\psi$ &  1.99594 &  1.99294 &  1.99083 & 1.98507 &  1.9919 &  1.99305 &  1.99307 &  1.99307\\
			\hline
			$\text{Err}_\mathbf{A}$ & 1.9982& 1.999& 2.07292 & 2.55124 & 1.36199 &  1.05806 &  1.01602 &  1.0111\\
			$\text{Err}_\gamma$& 2.00008 &1.99683 & 1.99408 &   1.98795 &  1.99617 & 1.99733 &  1.99735  & 1.99735\\
			\hline
			$\text{Err}_{\operatorname{\textbf{curl}}\gamma}$ &  2.00681  &2.00412 & 2.00279 &   2.00014 & 2.00375 &  2.00425 &  2.00425 &  2.00425\\
			\hline
			$\text{Err}_{\text{div}(\mathbf{A})}$  & 1.99169 &  2.00495 &  2.00072 &1.67622 &  0.941629 & 0.116245 &  -0.00869583 &  -0.0230222\\
			\hline
		\end{tabular}
	}
	\caption{2D manufactured solutions benchmark. Finite elements of order r = 1. Computed orders of $\psi$, $\mathbf{A}$, $\gamma$, $\operatorname{\textbf{curl}}\gamma$ and $\text{div}(\mathbf{A})$ for different gauges.}
	\label{tableDeg1allOm}
\end{table}

\begin{table}[h!]
	\hspace{-6em}
	\footnotesize
	{
		\begin{tabular}{lclclclclclcl}
			\hline
			&  $\omega = 1$ & $\omega = 10^{-1}$ & $\omega = 10^{-2}$ & $\omega = 10^{-3}$ & $\omega = 10^{-4}$ & $\omega = 10^{-5}$ &$\omega = 10^{-6}$ & $\omega = 0$\\ 
			\hline
			Richardson method& 1.9982& 1.999& 2.0729 & 2.5512 & 1.3619 &  1.0581 &  1.0160 &  1.0111\\
			\hline
			Graphic method& 1.9893 & 1.9897 & 2.1514 & 2.7523 & 1.6501 & 1.0781 & 1.0063 &  0.9982\\
			\hline
		\end{tabular}
	}
	\caption{2D manufactured solutions benchmark. Finite elements of order r = 1. Comparison  between the Richardson method and the graphical method for the estimation of convergence orders.}
	\label{tableDeg1VectPotCompGraphVsRich}
\end{table}

\begin{table}[h!]
	\hspace{-6em}
	\footnotesize
	{
		\begin{tabular}{lclclclclclcl}
			\hline
			&  $\omega = 1$ & $\omega = 10^{-1}$ & $\omega = 10^{-2}$ & $\omega = 10^{-3}$ & $\omega = 10^{-4}$ & $\omega = 10^{-5}$ &$\omega = 10^{-6}$ & $\omega = 0$\\ 
			\hline
			Richardson method & 1.9917 &  2.0049 &  2.0007 &1.6762 &  0.9416 & 0.1162 &  -0.0087&  -0.0230\\
			\hline
			Graphic method&  1.9790& 1.9903 & 1.9523 &  1.7748& 0.5894 & 0.01583 & -0.0635& -0.0643\\
			\hline
		\end{tabular}
	}
	\caption{2D manufactured solutions benchmark. Same caption as Tab. \ref{tableDeg1VectPotCompGraphVsRich} for $\operatorname{div}\mathbf{A}$.}
	\label{tableDeg1DivACompGraphVsRich}
\end{table}

In conclusion, the convergence orders are optimal when $\omega \geq 10^{-2}$. In the next section, we compute convergence orders for the case $r = 2$.

\subsection{Convergence analysis of the scheme \eqref{rtDiscrete} for the case $r = 2$}

Table \ref{table:3} shows the orders obtained with the Richardson extrapolation technique for the case $\omega = 1$. The space step is $\displaystyle \Delta x = \frac{1}{16} = 0.0625$, the time step is $\Delta t = 10^{-3}$ and we compute $N_{iter} = 125$ iterations. The results agree with \cite{gaoSun2015} except for the magnetic field $\gamma$; we observe an order equal to 4. We report on Tab. \ref{table:4} the results for different values of $\omega$. Like the case $r =1$, the quantities $\psi$, $\gamma$ and $\operatorname{\textbf{curl}} \gamma$ are not affected by the gauge choice. $\mathbf{A}$ (resp. $\operatorname{div}\mathbf{A}$) is losing one order (respectively two orders) when we decrease $\omega$ from 1 to 0.

\begin{table}[h!]
	\begin{tabular}{lclclclclclcl}
		\hline
		&  $\text{Err}_\psi$ &$\text{Err}_\mathbf{A}$ &$\text{Err}_\gamma$& $\text{Err}_{\operatorname{\textbf{curl}}\gamma}$ &$\text{Err}_{\text{div}(\mathbf{A})}$ \\
		\hline
		$||u_{\frac{\Delta x}{2}} - u_{\Delta x}||_{L^2} $ & $ 3.77341\cdot10^{-5}$ & $ 2.61501\cdot10^{-5}$ & $1.80296\cdot10^{-6}$ & $0.000212573$ & $0.000127969$\\
		\hline
		$||u_{\frac{\Delta x}{4}} -u_{\frac{\Delta x}{2}}||_{L^2}$&  $4.73961\cdot10^{-6}$ & $3.26017\cdot10^{-6}$ & $1.12463\cdot10^{-7}$ &$ 2.65818\cdot10^{-5}$ &  $ 1.60875\cdot10^{-5}$\\
		\hline
		$\text{order} $ & $2.99303$ & $3.0038$ & $4.00284$ &$2.99945$ & $2.99178$ \\
		\hline
	\end{tabular}
	\caption{2D manufactured solutions benchmark. Finite elements of order r = 2. Computed $L^2$ norm errors of $\psi$, $\mathbf{A}$, $\gamma$, $\operatorname{\textbf{curl}}\gamma$ and $\text{div}(\mathbf{A})$ for the Lorenz gauge.}
	\label{table:3}
\end{table}

\begin{table}[h!]
	\centering
	\footnotesize
	{
		\hspace{-6em}
		\begin{tabular}{lclclclclclcl}
			\hline
			&  $\omega = 1$ & $\omega = 10^{-1}$ & $\omega = 10^{-2}$ & $\omega = 10^{-3}$ & $\omega = 10^{-4}$ & $\omega = 10^{-5}$ &$\omega = 10^{-6}$ & $\omega = 0$\\ 
			\hline
			$\text{Err}_\psi$ &  2.99303 &  2.99293   &   2.99294 &  2.99311 &  2.993 &  2.99299 & 2.99295 &   2.99295 \\
			\hline
			$\text{Err}_\mathbf{A}$ & 3.0038 & 3.00705 &   3.22632 &  3.65457 &  2.32248& 2.0462 &  2.07952 &   1.98825 \\
			$\text{Err}_\gamma$&  4.00284 &  4.002  & 3.99905 &  3.98761 & 3.99725 & 4.0023 & 4.00289 &    4.00296\\
			\hline
			$\text{Err}_{\operatorname{\textbf{curl}}\gamma}$ &  2.99945  & 3 & 2.99963 &  2.99651  &  2.99804 &  2.99994 
			&   2.99993  &  2.99995 \\
			\hline
			$\text{Err}_{\text{div}(\mathbf{A})}$  & 2.99178 & 2.97989  &   2.98608 &   2.68217 & 1.31567 & 1.05971 &  1.22292 &  0.984674 \\
			\hline
		\end{tabular}
	}
	\caption{2D manufactured solutions benchmark. Finite elements of order r = 2. Computed orders of $\psi$, $\mathbf{A}$, $\gamma$, $\operatorname{\textbf{curl}}\gamma$ and $\text{div}(\mathbf{A})$ for different gauges.}
	\label{table:4}
\end{table}

In conclusion, as for the case $r = 1$, the convergence orders are optimal when $\omega \geq 10^{-2}$. This is consistent with the faster energy decrease displayed  in  Fig. \ref{genGfig4} for $\omega \geq 10^{-2}$. 

\section{Results in three dimensions of space}

In this section, we consider the mixed FE scheme \eqref{rtDiscrete3D} in three dimensions. We compute orders for different values of $\omega$ using manufactured solutions on the unit cube. Then, we study three domain configurations: the unit cube, a sphere and a sphere with a geometrical defect.

\subsection{Computation of convergence orders for the scheme \eqref{rtDiscrete3D}}

We consider  the following manufactured system for the TDGL model on the unit cube $(0,1)^3$:
\begin{plain}
	\begin{equation}
		\eqalign{
		&	\frac{\partial \psi}{\partial t}  -i\kappa\omega\text{div}(\mathbf{A}) 
		-\left( \frac{1}{\kappa}\nabla -i\mathbf{A}\right)^2\psi - \psi +|\psi|^2\psi = g,\cr
		&	 \frac{\partial {\mathbf A}}{\partial t} - \omega \nabla \text{div}(\mathbf{A}) + \operatorname{\textbf{curl}} \operatorname{\textbf{curl}} \mathbf{A}
		- \frac{1}{2i\kappa}\left( \psi^*\nabla\psi - \psi\nabla\psi^* \right)	+ |\psi|^2{\mathbf A} = \operatorname{\textbf{curl}} {\mathbf H} + \mathbf{f},
}
	\label{TDGL3DGeneralizedGArtificial}
\end{equation}
\end{plain}
with boundary and initial conditions
\begin{equation}
	\begin{array}{ll}
		\displaystyle	\nabla\psi \cdot \mathbf{n} = 0, \quad \operatorname{\textbf{curl}} \mathbf{A} \times\mathbf{n}= \mathbf{H}\times \mathbf{n}, \quad \omega\mathbf{A}\cdot \mathbf{n} = 0, \quad \mbox{on } \partial\Omega,
	\end{array}
	\label{TDGL3DBcGeneralizedGArtificial1}
\end{equation}	
where $\mathbf{f}$ and $g$ are defined such that the exact solution of \eqref{TDGLGeneralizedGArtificial} is
\begin{plain}
	\begin{equation}
		\eqalign{
	&	\psi = \exp(t)\cos(\pi y)\cos(\pi z) + i\exp(t)\cos(\pi x)\cos(\pi z), \cr
	&	\mathbf{A} = \left(
		\begin{array}{ll}
			\exp(t)\sin(\pi x) \sin(\pi y)\\
			\exp(t)\sin(\pi y)\sin(\pi z)\\
			\exp(t)\sin(\pi z)
		\end{array}\right),\cr
	&	\mathbf{H} = \left(
		\begin{array}{ll}
			-\pi	\exp(t)\sin(\pi y) \cos(\pi z)\\
			0\\
			-\pi	\exp(t)\sin(\pi x)\cos(\pi y)
		\end{array}\right).
}
\end{equation}
\end{plain}

Table \ref{table:5} shows the orders obtained with the Richardson extrapolation technique for $\displaystyle \Delta x = \frac{1}{10} = 0.1$. The time step is $\Delta t = 10^{-3}$ and we make $N_{iter} =100$ iterations. The results for the Lorenz gauge are in agreement with \cite{gaoSun2015} except for $\psi$. We observe an order 2 for the order parameter. Only $\operatorname{div}\mathbf{A}$ is affected by the gauge. The convergence rate of $\operatorname{div}\mathbf{A}$ increases for $\omega = 10^{-2}$ then decreases towards 0.

\begin{table}[h!]
	\hspace{-6em}
	\footnotesize
	{
		\begin{tabular}{lclclclclclclclclclcl}
			\hline
			&  $\omega = 1$ & $\omega = 10^{-1}$ & $\omega = 10^{-2}$ & $\omega = 10^{-3}$ & $\omega = 10^{-4}$ & $\omega = 10^{-5}$ &$\omega = 10^{-6}$  & $\omega = 0$\\ 
			\hline
			$\text{Err}_\psi$ & 1.95673 &    1.95367 &   1.95264 &  1.95284 & 1.9536 &  1.95362 & 1.95362 &  1.95362\\
			\hline
			$\text{Err}_\mathbf{A}$ & 1.00174 &  1.00173 &  0.999845 &  1.01192 & 1.00799 &  1.00236 & 1.00153 &  1.00143\\
			$\text{Err}_\gamma$&  0.995947 &   0.995988 & 0.997945 &   0.995992 &  0.995994 &  0.995994 & 0.995994 &  0.995994\\
			\hline
			$\text{Err}_{\operatorname{\textbf{curl}}\gamma}$ & 0.998499 &   0.998555 &  0.998559 &  0.995992 &  0.998569 &  0.99856 & 0.998569 & 0.998569\\
			\hline
			$\text{Err}_{\text{div}(\mathbf{A})}$  & 0.998665 &  1.00779 &  1.35479 &  0.901413 &  0.198305 & -0.00942684 &   -0.0358771 &   -0.0388907\\
			\hline
		\end{tabular}
	}
	\caption{3D manufactured solutions benchmark. Computed orders of $\psi$, $\mathbf{A}$, $\gamma$, $\operatorname{\textbf{curl}}\gamma$ and $\text{div}(\mathbf{A})$ for different gauges.}
	\label{table:5}
\end{table}

\vspace{1em}
In conclusion, as in the 2D case, the convergence orders are optimal when $\omega \geq 10^{-2}$. 

\subsection{Numerical examples for 3D configurations}

We use the scheme \eqref{rtDiscrete3D} to analyse the convergence with respect to the choice of the gauge in 3 cases: the unit cube, a sphere and a sphere with a geometrical defect.

\begin{itemize}
	\item The unit cube: we set $\kappa = 10$, $\mathbf{H} = (0,0,5)$ and $\Delta t = 0.1$. The mesh is uniform and the number of nodes per $\xi$ is 3. Figure \ref{cubeCompareGauge1} shows the vortex pattern at $t = 100$ for $\omega = 1$. For other values of $\omega$, the final state is identical. To highlight the difference between the different gauges, the ratios $\displaystyle \frac{|{\cal G}_{n+1} - {\cal G}_n|}{{\cal G}_n}$, $n = 0\cdots 100$, are plotted on Fig. \ref{cubeCompareGauge1}. We observe that the convergence is similar for $\omega = 1, 10^{-1}, 10^{-2}$. For lower values of $\omega$, a change of regime appears and the convergence is much slower.

	\begin{figure}[!h]
		\includegraphics[width = 0.4\textwidth]{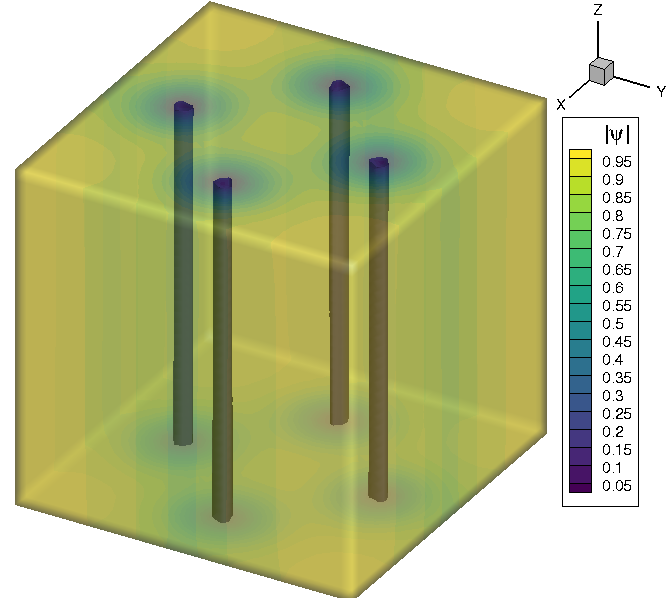}\hfill
		\includegraphics[width = 0.5\textwidth]{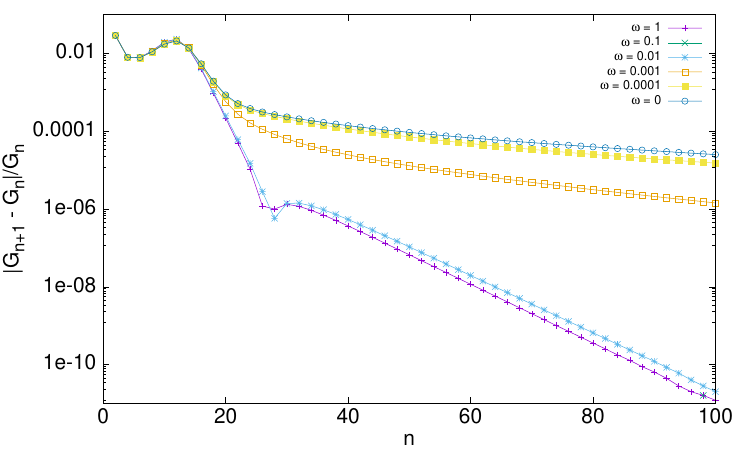}
		\caption{3D benchmark with unit cube domain. Density $|\psi|$ for $\kappa = 10$, $H = 5$, at $t = 100$ with $\omega = 1$ (left). Relative energy  differences $\displaystyle |{\cal G}_{n+1} - {\cal G}_n|/{\cal G}_n$ for different gauges.}
		\label{cubeCompareGauge1}
	\end{figure}
	\begin{figure}[!h]
	\includegraphics[width = 0.4\textwidth]{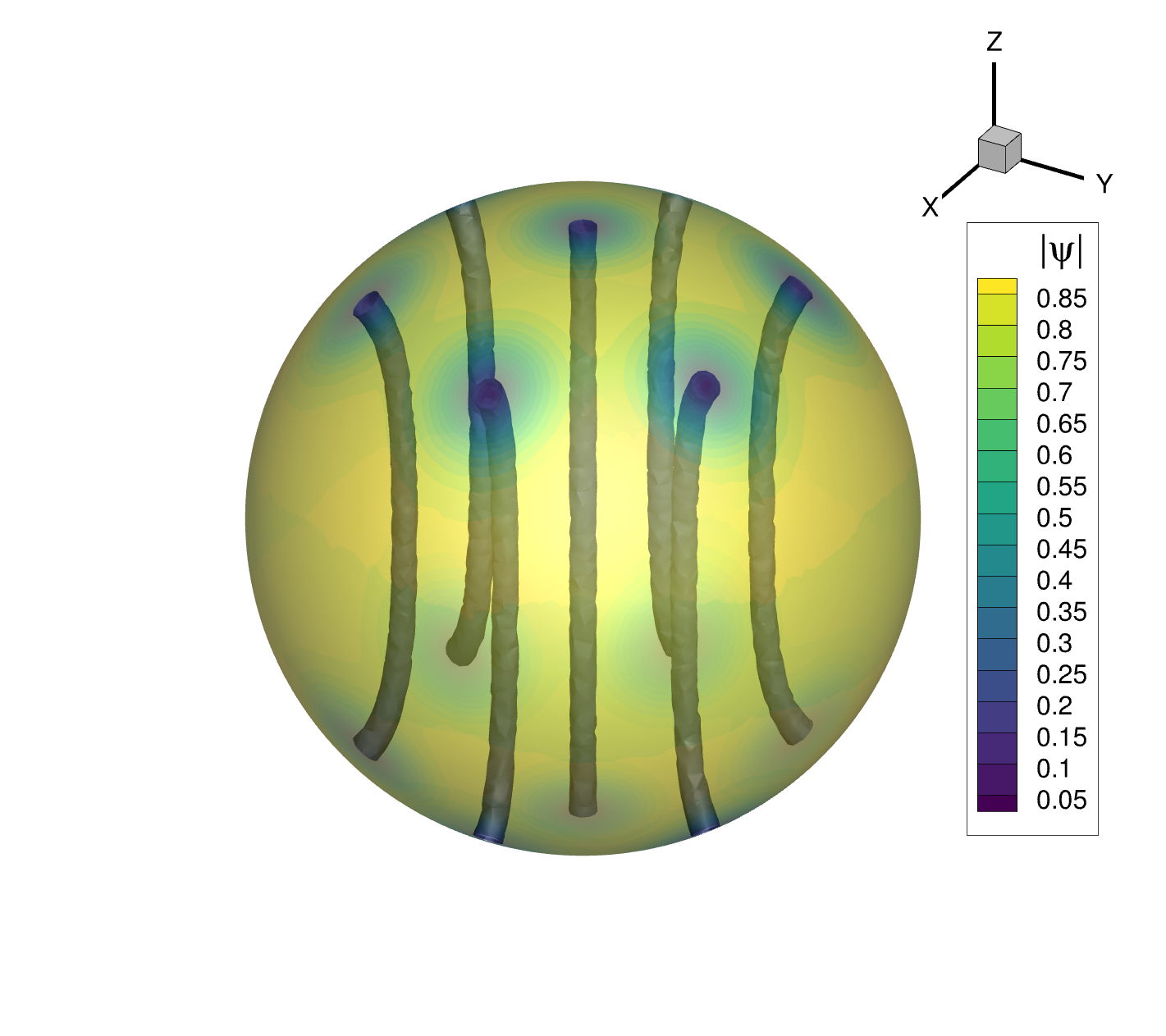}\hfill
	\includegraphics[width = 0.5\textwidth]{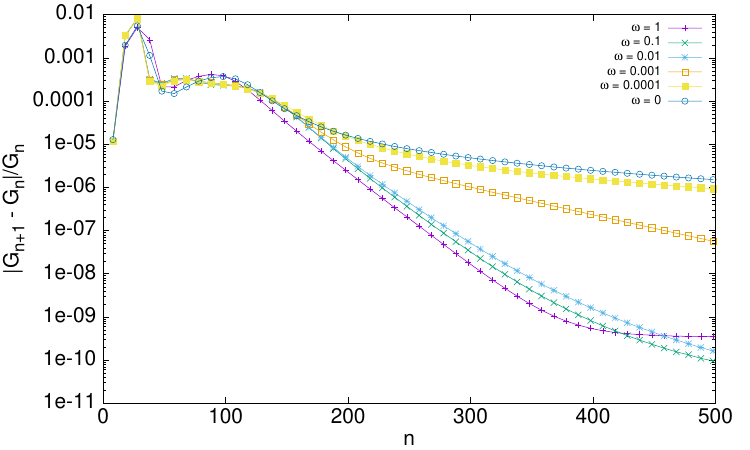}
	\caption{3D benchmark with spherical domain. Contours of the density $|\psi|$ for $\kappa = 10$, $H = 5$, at $t = 500$ with $\omega = 1$ (left). Relative energy  differences $\displaystyle |{\cal G}_{n+1} - {\cal G}_n|/{\cal G}_n$ for different gauges (right).}
	\label{sphereCompareGauge1}
\end{figure}

	\item A sphere with or without a geometrical defect:
	the domain is the sphere of radius $\displaystyle \frac{\sqrt{2}}{2}$ with or without a defect. We set $\kappa = 10$, $\mathbf{H} = (0,0,5)$ and $\Delta t = 0.1$. The mesh is uniform and the number of nodes per $\xi$ is 3. Figure  \ref{mesh3DPackB} shows the mesh for the sphere with a defect. Vortex patterns at equilibrium for $\omega = 1$ are shown on Figs.  \ref{sphereCompareGauge1} and \ref{packmanBCompareGauge1}. The patterns are identical for other choices of $\omega$. The ratios $\displaystyle \frac{|{\cal G}_{n+1} - {\cal G}_n|}{{\cal G}_n}$, $n = 0\cdots 500$, are plotted in Figs. \ref{sphereCompareGauge1} and \ref{packmanBCompareGauge1}. We observe a better convergence for the cases $\omega = 1, 10^{-1}, 10^{-2}$. These results are consistent with the convergence rates obtained from Tab. \ref{table:5}.

	\begin{figure}[!h]
		\centering
		\includegraphics[width = 0.4\textwidth]{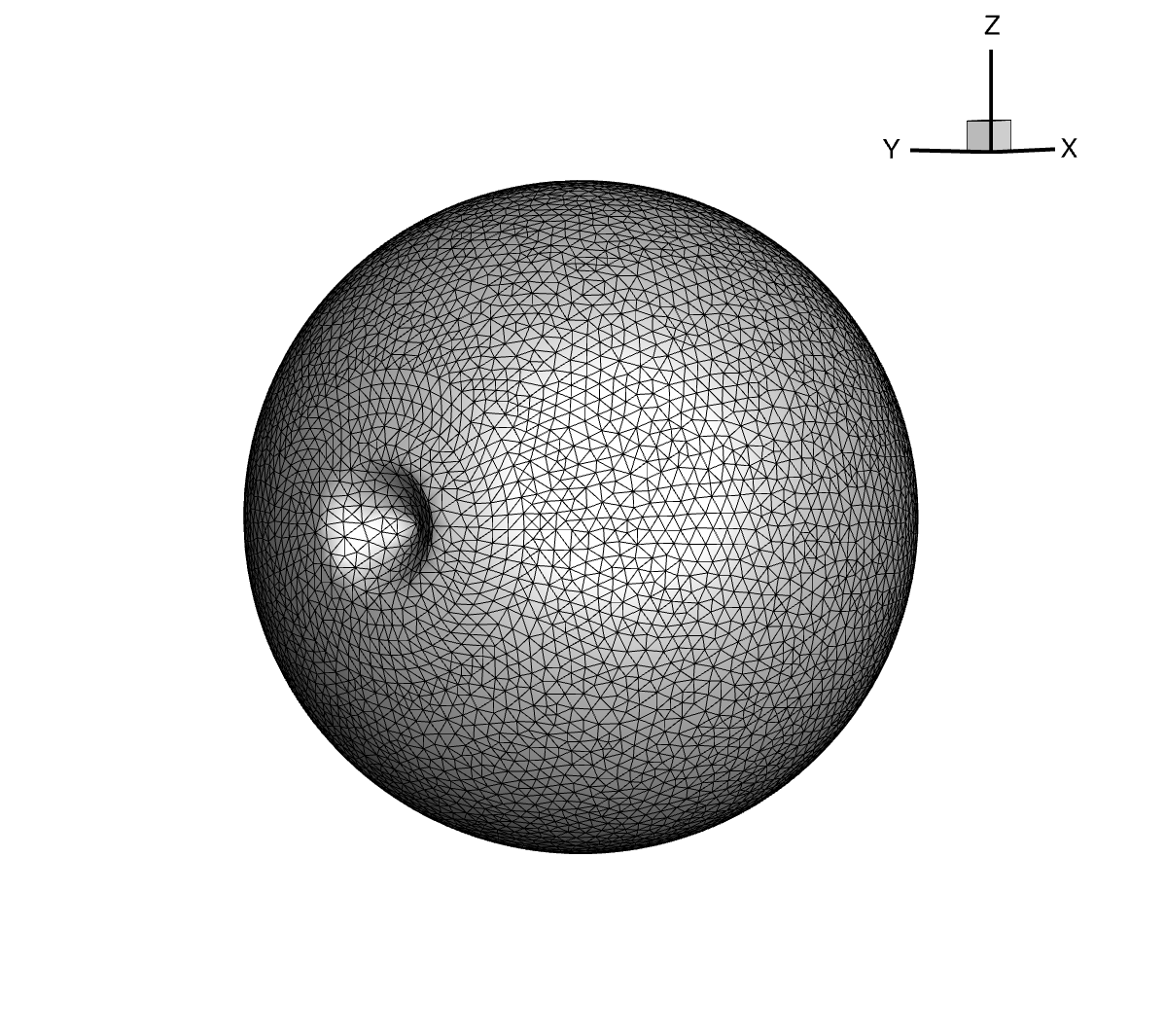}
		\caption{3D mesh of the sphere of radius $\frac{\sqrt{2}}{2}$ with a geometrical defect.}
		\label{mesh3DPackB}
	\end{figure}
	\begin{figure}[!h]
		\includegraphics[width = 0.4\textwidth]{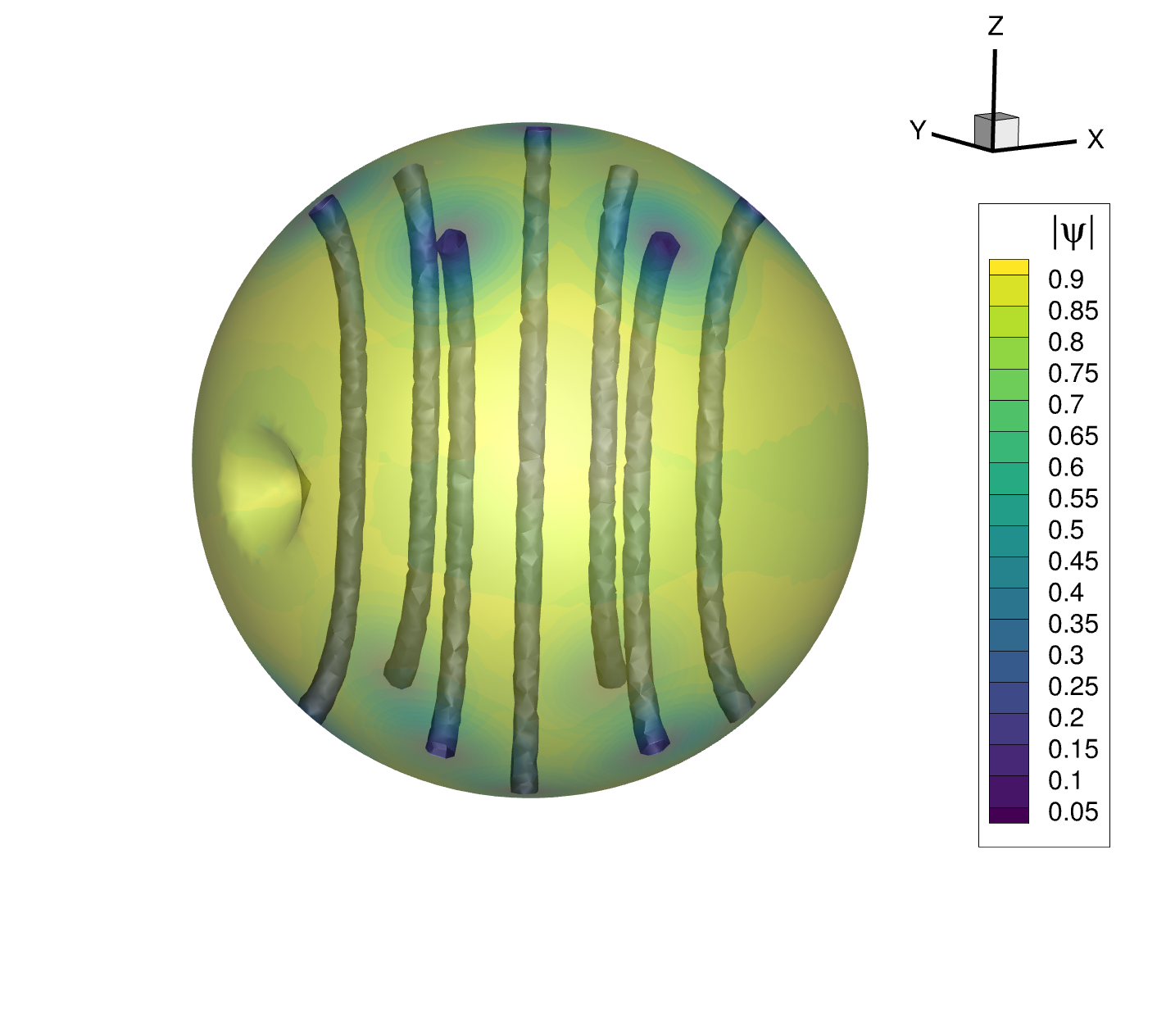}\hfill
		\includegraphics[width = 0.5\textwidth]{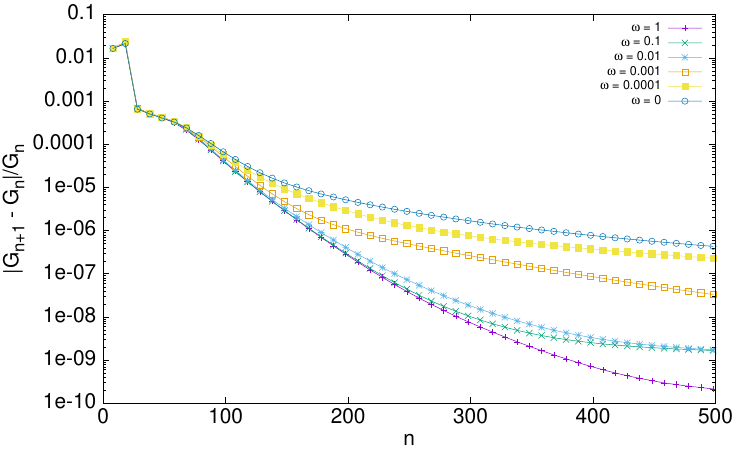}
		\caption{3D benchmark with spherical domain and a geometrical defect. Density $|\psi|$ for $\kappa = 10$, $H = 5$, at $t = 500$ with $\omega = 1$ (left). Relative energy  differences $\displaystyle |{\cal G}_{n+1} - {\cal G}_n|/{\cal G}_n$ for different gauges (right).}
		\label{packmanBCompareGauge1}
	\end{figure}
\end{itemize}

\section{Summary and conclusions}\label{sec-conclusions}

We presented a comparative study of different gauges for the TDGL model in the unified framework of the $\omega$-gauge theoretically introduced in \cite{fleckinger97}. Classical gauges were recovered from this model: $\omega = 1$ corresponds to the Lorenz gauge and $\omega = 0$ to the temporal gauge. We used the mixed finite element scheme introduced by \cite{gaoSun2015}, rewritten in the $\omega$-gauge. The present contribution is a first attempt, to the best of our knowledge, to numerically analyse the $\omega$-gauge formulation in FE settings.

First we studied a classical benchmark, the disk with a geometrical defect. For FE of order $r = 1$ and low values of the gauge parameter $\omega$, we observed a growing normal zone near the defect. This numerical artefact recalls the extended vortices found in \cite{sorensen10}. However, in our case, a finer resolution of the mesh or an increase of the FE order $r$ can solve the issue. As an alternative strategy to avoid such undesirable effects, we showed that using a higher value for $\omega$, typically above $10^{-2}$ is also efficient. This suggestion is confirmed  by plotting the energy decrease during computations, which is faster when $\omega \ge 10^{-2}$. Incidentally, by varying $\omega$, we also found a state of new lowest energy state (non-symmetrical with respect to the x-axis), which is similar to a minimizer of the renormalized energy introduced in \cite{sandier2008vortices} for a system of point vortices.

In the second part of our study, our goal was to assess the influence of $\omega$ on the convergence orders for $\psi$, $\mathbf{A}$, ${\gamma}$, $\operatorname{\textbf{curl}}{\gamma}$, and $\operatorname{div}\mathbf{A}$.
In 2D, only $\mathbf{A}$ and $\operatorname{div}\mathbf{A}$ were affected by a change of gauge. The degeneracy of the convergence orders, already observed in \cite{gao2016}, were recovered. In addition, we saw that the tipping point between optimal convergence and degeneracy occurred for $\omega$ between $10^{-2}$ and $10^{-3}$. Moreover, a careful study of the convergence curves on Fig. \eqref{orderGeneralizedG} showed that this tipping point also depends on the size of the mesh. These results are in agreement with the analysis of the previous 2D benchmark: increasing the gauge parameter $\omega$ or refining the mesh are the two ways to ensure the best convergence.

In 3D, the analysis was conducted for $r=0$. It appeared that only $\operatorname{div} \mathbf{A}$ was affected by a gauge choice. Quantities $\psi$, $\mathbf{A}$, $\boldsymbol{\gamma}$ and $\operatorname{\textbf{curl}}\boldsymbol{\gamma}$ were unaffected. The convergence rate was $1$, which is the theoretical value, except for $\psi$ for which we observed superconvegence with a rate equal to $2$. As in the 2D case, the degeneracy of convergence orders appeared for $\omega \le 10^{-3}$. 
Two new benchmarks, the sphere with and without a defect were analysed. Each benchmark showed a clear threshold value for $\omega$ between $10^{-3}$ and $10^{-2}$ consistent with our convergence analysis.

In conclusion, we analysed in detail the link between the choice of the gauge and the behaviour (convergence order) of mixed FE schemes used to solve the TDGL system of equations. A potential user of FE methods must be aware that numerical artefacts could appear in numerical simulations. We suggested several strategies to avoid such undesirable effects and tested them on 2D and 3D benchmarks. This suggests that 
configurations with geometries relevant for actual superconductors can be successfully simulated with the $\omega$-gauge formulation.

\section*{Statements and Declarations}

The authors declare no competing interests.

%


\end{document}